\documentclass{article}
\usepackage{latexsym}
  \newpage

\mathchardef\tnode="020E 

\def\arc{
  \hbox{\kern -0.15em
  \vbox{\hrule width 3.4em height 0.6ex depth -0.5 ex}
  \kern -0.33em}}

\def\darc{
  \rlap{\lower0.2ex\arc}{\raise0.2ex\arc}}

\def\stroke#1{
  \kern 0.05em
  \rlap\arc{{\textstyle{#1}}\atop\phantom\arc}
  \kern -0.22em}

\def\dstroke#1{
  \kern 0.05em
  \rlap\darc{{\textstyle{#1}}\atop\phantom\darc}
  \kern -0.22em}

\def\centerscript#1{
  \setbox0=\hbox{$\tnode$}
  \hbox to \wd0{\hss$\scriptstyle{#1}$\hss}}

\def\node{
  \def\super{}
  \def\sub{}
  \futurelet\next\dolabellednode}

  \let\sp=^
  \let\sb=_

  \def\dolabellednode{%
    \ifx\next\sb\let\next\getsub
    \else
      \ifx\next\sp\let\next\getsuper
      \else\let\next\donode
      \fi
    \fi
    \next}

  \def\getsub_#1{\def\sub{#1}\futurelet\next\dolabellednode}
  \def\getsuper^#1{\def\super{#1}\futurelet\next\dolabellednode}

  \def\donode{%
    \rlap{$\mathop{\phantom\tnode}\limits_{\centerscript{\sub}}
    ^{\centerscript{\super}}$}\tnode}

\def\varcdn{
  \kern -0.03em\vbox{\kern -0.5ex
  \hbox to \wd0{\hss\vrule width 0.04em depth 5.8ex\hss}
  \kern -0.3ex
  \hbox{$\tnode$}}}

\mathchardef\tnodef="020F 

\def\nodef{
  \def\super{}
  \def\sub{}
  \futurelet\next\dolabellednodef}

  \let\sp=^
  \let\sb=_

  \def\dolabellednodef{%
    \ifx\next\sb\let\next\getsubf
    \else
      \ifx\next\sp\let\next\getsuperf
      \else\let\next\donodef
      \fi
    \fi
    \next}

  \def\getsubf_#1{\def\sub{#1}\futurelet\next\dolabellednodef}
  \def\getsuperf^#1{\def\super{#1}\futurelet\next\dolabellednodef}

  \def\donodef{%
    \rlap{$\mathop{\phantom\tnodef}\limits_{\centerscript{\sub}}
    ^{\centerscript{\super}}$}\tnodef}

\def\varcdnf{
  \kern -0.03em\vbox{\kern -0.5ex
  \hbox to \wd0{\hss\vrule width 0.04em depth 5.8ex\hss}
  \kern -0.3ex
  \hbox{$\tnodef$}}}

\newtheorem{main}{Theorem}
\newtheorem{thm}{Theorem}[section]
\newtheorem{prop}[thm]{Proposition}
\newtheorem{defi}[thm]{Definition}
\newtheorem{lem}[thm]{Lemma}
\newtheorem{cor}[thm]{Corollary}

\newcommand\pf{\noindent{\bf Proof.}~~} \newcommand\qed{~~\hfill$\Box$}

\newcommand\Gm{\Gamma}
\newcommand\Dl{\Delta}
\newcommand\Lm{\Lambda}

\newcommand\Sg{\Sigma}
\newcommand\Th{\Theta}

\newcommand\cA{{\cal A}}

\newcommand\cC{{\cal C}}
\newcommand\cD{{\cal D}}
\newcommand\cE{{\cal E}}
\newcommand\cF{{\cal F}}
\newcommand\cH{{\cal H}}
\newcommand\cG{{\cal G}}
\newcommand\cP{{\cal P}}
\newcommand\cQ{{\cal Q}}
\newcommand\cR{{\cal R}}
\newcommand\cS{{\cal S}}
\newcommand\cU{{\cal U}}
\newcommand\cY{{\cal Y}}
\newcommand\la{{\langle}}
\newcommand\ra{{\rangle}}

\newcommand\tE{^2\!E_6(2)}

\title{Extended $F_4$-buildings and the Baby Monster}
\author{A.A. Ivanov \thanks{A part of the research was conducted
whilst this author was visiting Department of Mathematics and
Statistics, Bowling Green State University.}, D.V. Pasechnik and S.V. Shpectorov}

\begin{document}
\maketitle

\begin{abstract}

Let $\Th$ be the Baby Monster graph which is the graph on the set of
$\{3,4\}$-transpositions in the Baby Monster group $B$ in which two
such transpositions are adjacent if their product is a central
involution in $B$. Then $\Th$ is locally the commuting graph of central
(root) involutions in $\tE$. The graph $\Th$ contains a family of
cliques of size 120. With respect to the incidence relation defined via
inclusion these cliques and the non-empty
intersections of two or more of them form a geometry $\cE(B)$ with
diagram
$$c.F_4(t):~~~\node_1\stroke{\rm c}\node_2\arc\node_2\darc\node_t\arc\node_t$$
for $t=4$ and the action of $B$ on $\cE(B)$ is flag-transitive. We show
that $\cE(B)$ contains subgeometries $\cE($$^2E_6(2))$ and
$\cE(Fi_{22})$ with diagrams $c.F_4(2)$ and $c.F_4(1)$.
The stabilizers in $B$ of these subgeometries induce on them
flag-transitive actions of $^2E_6(2):2$ and $Fi_{22}:2$, respectively.
The geometries $\cE(B)$, $\cE($$^2E_6(2))$ and $\cE(Fi_{22})$ possess
the following properties: $(a)$ any two elements of type 1 are incident to at most 
one common element of
type 2 and $(b)$ three elements of type 1 are pairwise incident to common
elements of type 2 if and only if they are incident to a common element
of type 5. The paper addresses the classification problem of
$c.F_4(t)$-geometries satisfying $(a)$ and $(b)$. We construct three further examples
for $t=2$ with flag-transitive automorphism groups isomorphic to $3 \cdot
$$^2E_2(2):2$, $E_6(2):2$ and $2$$^{26}.F_4(2)$ and one for
$t=1$ with flag-transitive automorphism group $3 \cdot Fi_{22}:2$.
We also study the graph of an arbitrary
(non-necessary flag-transitive) $c.F_4(t)$-geometry satisfying $(a)$ and
$(b)$ and obtain a complete list of possibilities for the isomorphism
type of subgraph induced by the common neighbours of a pair of
vertices at distance 2. Finally, we prove that $\cE(B)$ is the only
$c.F_4(4)$-geometry, satisfying $(a)$ and $(b)$.

\end{abstract}

\section{Introduction} \label{s1}

The paper contributes to the geometric theory of sporadic simple
groups. Our notation and terminology is mostly standard (see
\cite{Pasi94} and \cite{I99}). On a diagram the types of elements increase
rightwards from 1 to the rank of the geometry. If an element of type 2
in a geometry $\cG$ is incident to exactly two elements of type 1 then
in the {\it graph} $\Gm(\cG)$ of $\cG$ the vertices and edges are the
elements of type 1 and 2 in $\cG$ subject to the natural incidence
relation (in general $\Gm(\cG)$ might contain multiple edges). Recall that if $\Gm$ is
a graph and $x\in \Gm$ is a vertex, then $\Gm_i(x)$ is the set
of vertices at distance $i$ from $x$ in $\Gm$. Sometimes we write
$\Gm(x)$ instead of $\Gm_1(x)$. If $\Sigma$ is a subset of the vertex-set of
$\Gm$, then $\Sigma$ also denotes the subgraph in $\Gm$ induced by
$\Sigma$. If for every $x\in \Gm$ the subgraph in
$\Gm(x)$ is isomorphic to a fixed graph $\Dl$ then $\Gm$ is said to be {\it
locally} $\Dl$. A graph $\Gm$ whose vertices are involutions in a 
group $G$ is said to be the {\it commuting graph} if $a,b \in \Gm$ are adjacent
if and only if the order $o(ab)$ of the product of $a$ and $b$ is 2
(equivalently if $a$ and $b$ commute). A graph $\Gm$ is said to be a
{\it $m$-clique extension} of a graph $\Dl$ if there is a mapping $\psi$ of
the vertex-set of $\Gm$ onto the vertex-set of $\Dl$ such that
$|\psi^{-1}(x)|=m$ for every $x \in \Dl$ and two distinct vertices $u,v$ in
$\Gm$ are adjacent if and only if the images $\psi(u)$ and $\psi(v)$
are either equal or adjacent in $\Dl$. A graph $\Gm$ is said to be the
{\it distance $1$-or-$2$ graph} of $\Dl$ if there is a bijection $\chi$ of
the vertex-set of $\Gm$ onto the vertex-set of $\Dl$ with the following
property: two vertices $x$ and $y$ in $\Gm$ are adjacent
if and only if the vertices
$\chi(x)$ and $\chi(y)$ are adjacent or at distance 2 in $\Dl$. By a
$2$-path in a graph $\Gm$ we mean an ordered sequence $(x,y,z)$ of
vertices, such that $x,z \in \Gm(y)$ and $z \in \Gm_2(x)$.
If $G$ is a group, acting by permutations on a set $X$, then for $x \in
X$ and $Y \subseteq X$ by $G(x)$ and $G[Y]$ we denote the stabilizer of
$x$ and (the setwise) stabilizer of $Y$ in $G$ (the upper and lower
indexes will have various different meanings). 

\vspace{1cm}

\noindent
\unitlength 0.80mm
\linethickness{0.4pt}
\begin{picture}(169.76,122.16)
\put(20.91,71.06){\circle{10.00}}
\put(115.52,106.95){\line(0,-1){73.33}}
\put(20.91,71.06){\makebox(0,0)[cc]{\small $a$}}
\put(63.74,80.50){\makebox(0,0)[cc]{\scriptsize 1782+44\,352}}
\put(58.76,59.78){\makebox(0,0)[cc]{\small $2^{2+20}.U_6(2):2$}}
\put(95.81,80.97){\makebox(0,0)[cc]{\scriptsize 2\,097\,152}}
\put(109.03,98.24){\makebox(0,0)[cc]{\scriptsize 3510}}
\put(114.73,122.16){\makebox(0,0)[cc]{\scriptsize 142\,155+694\,980}}
\put(125.75,100.36){\makebox(0,0)[cc]{\scriptsize 3\,127\,410}}
\put(123.98,48.64){\makebox(0,0)[cc]{\scriptsize 663\,552}}
\put(161.34,80.21){\makebox(0,0)[cc]{\scriptsize 69\,615}}
\put(157.68,57.10){\makebox(0,0)[cc]{\scriptsize 3\,898\,440}}
\put(135.64,39.00){\makebox(0,0)[cc]{\scriptsize 8064}}
\put(115.79,17.84){\makebox(0,0)[cc]
{\scriptsize 1+1134+36\,288+1\,161\,216+2\,097\,152}}
\put(95.41,60.14){\makebox(0,0)[cc]{\scriptsize 1\,824\,768}}
\put(106.38,43.73){\makebox(0,0)[cc]{\scriptsize 648}}
\put(144.76,25.50){\makebox(0,0)[cc]{\small $2^{1+20}.U_4(3):2^2$}}
\put(134.90,112.27){\makebox(0,0)[cc]{{\small $Fi_{22}:2$}}}
\put(160.89,91.12){\makebox(0,0)[cc]{{\small $F_4(2) \times 2^2$}}}
\put(77.74,75.74){\line(1,1){31.12}}
\put(121.31,33.73){\line(1,1){31.89}}
\put(108.86,33.34){\line(-1,1){31.51}}
\put(70.93,70.44){\oval(18.60,11.63)[]}
\put(25.97,70.83){\line(1,0){35.66}}
\put(114.72,112.68){\oval(18.60,11.63)[]}
\put(115.50,27.80){\oval(18.60,11.63)[]}
\put(160.46,70.83){\oval(18.60,11.63)[]}
\put(71.32,70.83){\makebox(0,0)[cc]{\small $\Th_1(a)$}}
\put(114.72,112.68){\makebox(0,0)[cc]{\small $\Th_2^3(a)$}}
\put(115.50,28.19){\makebox(0,0)[cc]{\small $\Th_2^4(a)$}}
\put(160.46,70.83){\makebox(0,0)[cc]{\small $\Th_3(a)$}}
\put(35.66,74.32){\makebox(0,0)[cc]{\scriptsize 3\,968\,055}}
\put(58.53,74.32){\makebox(0,0)[cc]{\scriptsize 1}}
\put(20.74,60.00){\makebox(0,0)[cc]{\small $2 \cdot $$^2E_6(2):2$}}
\put(115.56,7.41){\makebox(0,0)[cc]{Fig. 1}}
\end{picture}

\bigskip

The sporadic simple group $B$ known as the {\it Baby Monster} contains a
conjugacy class $\Th$ of involutions (with centralizers of the form $2
\cdot \tE:2$) such that for any two involutions
$a,b \in \Th$ the order $o(ab)$ of their product is at most 4. This is
to say, $\Th$ is a class of $\{3,4\}$-transpositions in $B$.
Let $\Th$ denote also the graph on this set, in which two transpositions are
adjacent whenever their product belongs to the class of central
involutions in $B$ (with centralizers of the form $2^{1+22}_+.Co_2$).
We call $\Th$ the {\em Baby Monster graph}. The suborbit diagram of
$\Th$ with respect to the action of $B$ is given on Fig 1.

\medskip

Let $a \in \Th$ and $B(a)=C_{B}(a)$ be the stabilizer of $a$ in the
action of $B$ on $\Th$ by conjugation. The diameter of $\Th$ is 3;
$B(a)$ acts transitively on $\Th_1(a)$; has two orbits $\Th_2^3(a)$ and
$\Th_2^4(a)$ on $\Th_2(a)$, where
$$\Th_2^m(a)=\{b \mid b \in \Th_2(a), o(ab)=m\};$$
$a$ commutes with every involution from $\Th_3(a)$ and $B(a)$ acts
transitively on this set. This implies particularly that $a$
(considered as an element of $B$) fixes $\{a\} \cup \Th_1(a) \cup
\Th_3(a)$ and acts on $\Th_2(a)$ with orbits of length 2.

\medskip

Locally $\Th$ is the commuting graph $\Dl$ of central (root) involutions
in $B(a)/\la a \ra \cong \tE:2$. In terms of the $F_4$-building
$\cF$ with the digram
$$\node_2\arc\node_2\darc\node_4\arc\node_4$$
the vertices of $\Dl$ are the elements of type 1
with two of them adjacent if incident to a common element of type 4.

The group $B$ contains a maximal 2-local subgroup $S \sim
2^{9+16}.Sp_8(2)$. The centre of $O_2(S)$ contains exactly 120
$\{3,4\}$-transpositions which induce in $\Th$ a clique (a maximal
complete subgraph) $Q$. Let $\cE(B)$ be an incidence system formed by
the set $\cQ$ of images of $Q$ under $B$ together with the
non-empty intersections of two or more cliques from $\cQ$; the
incidence relation is via inclusion and the type of an element is
determined by its size. Then $\cE(B)$ is a geometry
with the diagram
$$c.F_4(t):~~~\node_1\stroke{\rm c}\node_2\arc\node_2\darc\node_t\arc\node_t$$
for $t=4$ mentioned in \cite{B85} and $B$
induces on $\cE(B)$ a flag-transitive action. The elements of type 1
are the vertices of $\Th$ and the residue of such a
vertex $a$ is isomorphic to the building $\cF$.

\medskip

In \cite{I92} within the characterization proof for the rank 5 Petersen
geometry of $B$ the following result has been established.

\begin{prop} \label{a}
Let $\Gm$ be a graph which is locally the commuting graph of central
involutions in $\tE$. Suppose that $\Gm$ possesses an automorphism
group $G$ such that for every $x\in\Gm$ the stabilizer $G(x)$ of $x$ in
$G$ induces on $\Gm(x)$ an action containing $\tE$. Then $\Gm$
is isomorphic to the Baby Monster graph $\Th$.\qed
\end{prop}

The main goal of this paper is to prove the assertion in
Proposition~\ref{a} without assuming anything about the 
automorphism group of $\Gm$ 
(see Theorem~\ref{aaa} below).
One of the intended applications of Theorem~\ref{aaa} is using it as
an identification tool in the classification of (non-necessary
flag-transitive) $P$-geometries of rank 5,
as a continuation of the project started in \cite{HS00}.

\medskip

Let $a \in \Th$ and $b \in \Th_2^3(a)$. Then the maximal intersection
with $\Th_3(a)$ (resp. with $\Th_3(a) \cap \Th_3(b)$) of a clique from
$\cQ$ is of size 64 (resp. 36). Let $\cE(\tE)$ and $\cE(Fi_{22})$
be the incidence systems formed by the maximal intersections of the
sets $\Th_3(a)$ and $\Th_3(a) \cap \Th_3(b)$, respectively with the
cliques from $\cQ$ together with the elements of $\cE(B)$ contained in
these sets. The incidence relation and type function are as in
$\cE(B)$. 

\begin{prop} \label{b}
\begin{itemize}
\item[{\rm (i)}] $\cE(\tE)$ is a $c.F_4(2)$-geometry, 
on which $B(a)/\la a \ra \cong $$^2E_6(2):2$ induces a flag-transitive action;

\item[{\rm (ii)}] $\cE(Fi_{22})$ is a $c.F_4(1)$-geometry, on which $B(a) \cap B(b) \cong Fi_{22}:2$ induces a
flag-transitive action.
\end{itemize}
\end{prop}

By the construction $\cE(\tE)$ is a subgeometry in $\cE(B)$ and
$\cE(Fi_{22})$ is a subgeometry in $\cE(\tE)$. By (\ref{a}) the
geometry $\cE(B)$ is simply connected which is not the case for
$\cE(\tE)$ and $\cE(Fi_{22})$. Let $\cA$ be the amalgam of maximal
parabolics associated with the action of $\tE:2$ on $\cE(\tE)$. Then 
$\cA$ is embedded into the unique non-split extension $3
\cdot \tE:2$ and hence there is a (3-fold) covering
$$\varphi:\cE(3 \cdot \tE) \to \cE(\tE)$$
of geometries. The preimage of $\cE(Fi_{22})$ with respect to $\varphi$
is connected and it is a flag-transitive triple cover $\cE(3 \cdot Fi_{22})$ of
$\cE(Fi_{22})$.

The subdegrees of the action of $E_6(2)$ on the cosets of $F_4(2)$
calculated in \cite{L93} show that one of the orbitals of this action
is the graph of a flag-transitive $c.F_4(2)$-geometry $\cE(E_6(2))$.
In a similar way the orbit lengths of $F_4(2)$ on the vectors of a
26-dimensional $GF(2)$-module calculated in \cite{CC88} prove existence
of a $c.F_4(2)$-geometry $\cE(2^{26}.F_4(2))$.

\medskip

We are focussed on the $c$-extensions $\cE$ of buildings
with diagrams $c.F_4(t)$, satisfying the following two conditions:

\begin{itemize}
\item[{\rm $(a)$}] any two elements of type 1 are incident to at most one
common element of type 2;
\item[{\rm $(b)$}] three elements of type 1 are pairwise incident to
common elements of type 2 if and only if they are incident to a common
element of type 5.
\end{itemize}

Let $\Gm=\Gm(\cE)$ be the graph of $\cE$. Then the above conditions are
equivalent, respectively, to the following ones:

\begin{itemize}
\item[{\rm $(a^\prime)$}] $\Gm$ has no multiple edges;
\item[{\rm $(b^\prime)$}] for an element $z$ of type 5 the vertices in $\Gm$
incident to $z$ induce a complete graph $Q(z)$ and every triangle in
$\Gm$ is contained in $Q(z)$ for some $z$.
\end{itemize}

The following theorems constitute the main result of the paper.

\begin{main} \label{aaa}
Let $\Gm$ be a graph which is locally the commuting graph of the
central involutions in $^2E_6(2)$. Then $\Gm$ is isomorphic to the Baby
Monster graph $\Theta$.
\end{main}

\begin{main} \label{bbb}
Let $\cE$ be a $c.F_4(4)$-geometry satisfying the conditions (a) and
(b). Then $\cE$ is the Baby Monster geometry $\cE(B)$.
\end{main}

We emphasise again that in Theorems~\ref{aaa} and \ref{bbb} no
automorphisms are assumed {\it a priori}. The condition $(a)$ excludes
some possible degenerate cases and it holds in every flag-transitive
$c.F_4(t)$-geometry. On the other hand the examples presented at the
end of Section~\ref{s4} show that the class of $c.F_4(t)$-geometries
which do not satisfy $(b)$ is rather large and even in the
flag-transitive case it is difficult to expect the complete
classification. 

\medskip

The paper is organized as follows. In Section~\ref{s2} we recall some
basic properties of the $F_4$-buildings with 3 points per a line. In
Section~\ref{s3} we discuss some general properties of
$c.F_4(t)$-geometries and their graphs. In Section~\ref{s4} we present
and discuss some examples of $c.F_4(t)$-geometries. In
Section~\ref{s5} we construct an important class of subgraphs in a
graph of $c.F_4(t)$-geometry. In Sections~\ref{s6} and \ref{s7} we
characterize the $\mu$-graphs which are subgraphs induced by the
common neighbours of a pair of vertices at distance 2. Finally in
Section~\ref{s8} (which consists of a few subsections) we restrict
ourselves to the case $t=4$ and eventually prove the main theorems.

\section{Some $F_4$-buildings} \label{s2}
Let $\cF=$$^t\cF$ be a building with diagram
$$\node_2\arc\node_2\darc\node_t\arc\node_t$$
By the classification \cite{T74} of spherical buildings we obtain that
$t=4,2$ or $1$ and the (flag-transitive) automorphism group $F=$$^tF$ of
$\cF$ is isomorphic to
$$\tE:S_3,~~F_4(2) {\rm ~~~or~~~} \Omega_8^+(2):S_3,$$
respectively. The
elements of type 1, 2, 3 and 4 in $\cF$ will also be called {\it
points}, {\it lines}, {\it planes} and {\it symplecta}, respectively.
The points of $\cF$ are the central (root) involutions in $F$ for
$t=1$ and $4$ and one of the two classes of such involutions for $t=2$. Let
$^t\Psi$ be the collinearity graph of $\cF$, that is the graph on points
in which two of them are adjacent if incident to a common line. The
following result is well known (\cite{Coh83} and \cite{Coo83}).

\begin{lem} \label{1}
The suborbit diagram of $\Psi$ is the one given on Fig.~$2$,  where $\{p\}$,
$\Psi_1(p)$, $\Psi_2^2(p)$, $\Psi_2^4(p)$ and $\Psi_3(p)$ are the
orbits of $F(p):=C_F(p)$ on $\Psi$ ($p$ is a point). Furthermore
\begin{itemize}
\item[{\rm (i)}] $F(p) \cong 2^{1+20}_+:U_6(2).S_3$,
$2^{1+6+8}:Sp_6(2)$ and $2^{1+8}_+:(S_3 \wr S_3)$ for $t=4,2$ and $1$;
\item[{\rm (ii)}] $q \in \Psi_1(p)$ if and only if $q \in O_2(F(p))$;
\item[{\rm (iii)}] $q \in \Psi_2^2(p)$ if and only if $q \in F(p)
\setminus O_2(F(p))$, in which case there is a unique symplecton
incident to $p$ and $q$;
\item[{\rm (iv)}] $q \in \Psi_2^4(p)$ if and only if $o(pq)=4$, in
which case $(pq)^2$ is the unique point collinear to both $p$ and $q$;
\item[{\rm (v)}] $q \in \Psi_3(p)$ if and only if $o(pq)=3$, in which
case $F(p) \cap F(q)$ is a complement to $O_2(F(p))$ in $F(p)$, if
$\{p,r_1,r_2\}$ is a line containing $p$, then up to renumbering $r_1
\in \Psi_2^4(q)$, $r_2 \in \Psi_3(q)$;
\item[{\rm (vi)}] the orbits of $O_2(F(p))$ on $\Psi_1(p)$,
$\Psi_2^2(p)$, $\Psi_2^4(p)$ and $\Psi_3(p)$ are of length $2$,
$2^4\cdot t$, $2^5 \cdot t^3$ and $2^9 \cdot t^6=|\Psi_3(p)|$, respectively;
\item[{\rm (vii)}] the orbits of $O_2(F(p))$ on $\Psi_1(p)$ and
$\Psi_2^4(p)$ naturally correspond to the lines incident to $p$ while
the orbits of $O_2(F(p))$ on $\Psi_2^2(p)$ naturally correspond to the
symplecta incident to $p$.\qed
\end{itemize}
\end{lem}

The points lines, planes and symplecta in $\cF$ can be identified with
the sets (of size 1, 3, 7 and $28t+7$, respectively) of points incident
to these elements, so that the incidence relation is via inclusion. The
geometry $\cF$ can be reconstructed from $\Psi$ in the following way.
If $q \in \Psi_2^2(p)$, then the unique symplecton containing $p$ and
$q$ is the minimal geodetically closed subgraph in $\Psi$ containing
$p$ and $q$. The remaining elements of $\cF$ are the non-empty
intersections of two or more symplecta.

\medskip

There are the following embeddings among the three geometries:
$$^1\cF \subset ^2\!\cF \subset ^4\!\cF$$
which can be seen as follows. The group $^4F$ contains an (outer)
involution $x$ such that $C_{^4F}(x) \cong \la x \ra \times $$^2F$.
Furthermore there is a conjugate $y$ of $x$ in $^4F$ such that $\la x,y
\ra \cong S_3$ and $C_{^4F}(\la x,y \ra) \cong $$^1F$. The subgraph in
$^4\Psi$ induced by the vertices fixed by $x$ is isomorphic to $^2\Psi$
and the subgraph induced by the vertices fixed by $\la x,y \ra$
(equivalently by $xy$) is isomorphic to $^1\Psi$.

\medskip

\unitlength 0.85mm
\linethickness{0.4pt}
\begin{picture}(167.78,147.97)
\put(85.33,19.63){\circle{8.00}}
\put(84.66,48.63){\oval(39.33,8.00)[]}
\put(40.00,89.96){\oval(39.33,8.00)[]}
\put(85.33,139.30){\oval(30.00,8.00)[]}
\put(85.33,23.63){\line(0,1){21.00}}
\put(57.33,93.96){\line(2,5){16.40}}
\put(57.00,85.96){\line(1,-3){11.11}}
\put(101.66,52.63){\line(1,3){11.00}}
\put(59.66,89.96){\line(1,0){50.67}}
\put(85.33,19.63){\makebox(0,0)[cc]{\small 1}}
\put(85.33,48.63){\makebox(0,0)[cc]{\small $6(2t^2+1)(2t+1)$}}
\put(41.00,89.96){\makebox(0,0)[cc]{\small $96t^3(2t^2+1)(2t+1)$}}
\put(134.00,90.29){\makebox(0,0)[cc]{\small $16t(2t^2+1)(t^2+t+1)$}}
\put(85.33,139.63){\makebox(0,0)[cc]{\small $512t^6$}}
\put(100.33,28.30){\makebox(0,0)[cc]{\scriptsize $6(2t^2+1)(2t+1)$}}
\put(88.66,40.96){\makebox(0,0)[cc]{\scriptsize 1}}
\put(85.33,56.30){\makebox(0,0)[cc]{\scriptsize $4(t^2+t+1)+1$}}
\put(59.33,58.63){\makebox(0,0)[cc]{\scriptsize $16t^3$}}
\put(54.66,80.63){\makebox(0,0)[cc]{\scriptsize 1}}
\put(114.66,58.63){\makebox(0,0)[cc]{\scriptsize $8t(t^2+t+1)$}}
\put(119.66,80.63){\makebox(0,0)[cc]{\scriptsize $3(2t+1)$}}
\put(35.99,97.63){\makebox(0,0)[cc]{\scriptsize $(8t+2)(t^2+t+1)+1$}}
\put(133.66,97.30){\makebox(0,0)[cc]{\scriptsize $3(2t+1)$}}
\put(71.66,92.96){\makebox(0,0)[cc]{\scriptsize $2(t^2+t+1)$}}
\put(99.00,92.96){\makebox(0,0)[cc]{\scriptsize $12t^2(2t+1)$}}
\put(53.00,106.96){\makebox(0,0)[cc]{\scriptsize $16t^3$}}
\put(56.99,128.30){\makebox(0,0)[cc]{\scriptsize $3(2t^2+1)(2t+1)$}}
\put(85.00,147.97){\makebox(0,0)[cc]{\scriptsize $3(2t^2+1)(2t+1)$}}
\put(132.83,89.96){\oval(45.00,8.00)[]}
\put(96.33,19.96){\makebox(0,0)[cc]{\small $\{p\}$}}
\put(8.89,90.00){\makebox(0,0)[cc]{\small $\Psi_2^4(p)$}}
\put(167.78,90.00){\makebox(0,0)[cc]{\small $\Psi_2^2(p)$}}
\put(111.48,139.26){\makebox(0,0)[cc]{\small $\Psi_3(p)$}}
\put(115.67,48.30){\makebox(0,0)[cc]{\small $\Psi_1(p)$}}
\put(85.19,10.00){\makebox(0,0)[cc]{Fig. 2}}
\end{picture}

Let $\widehat \Psi$ be a graph on the point-set of $\cF$
in which two points $p$ and $q$ are adjacent if they are incident to a common
symplecton but not to a common line, equivalently if $q \in
\Psi_2^2(p)$. The suborbit diagram of $\widehat \Psi$ computed from that of
$\Psi$ is presented on Fig. 3.

\medskip

Let $\Dl=$$^t\Dl$ be the graph on the point-set of $\cF$ in which two
points are adjacent whenever they are incident to a common symplecton.
This means that two points $p$ and $q$ are adjacent in $\Dl$ if $q \in
\Psi_1(p) \cup \Psi_2^2(p)$, so that the edge-set of $\Dl$ is the union
of the edge-sets of $\Psi$ and $\widehat \Psi$. By (\ref{1}) $\Dl$ is the
commuting graph on points of $\cF$ (considered as involutions in $F$).

\medskip

The next two lemmas follow directly from the diagrams of $\Psi$ and
$\widehat \Psi$.

\begin{lem} \label{1.5}
If $(p,r,q)$ is a $2$-path in $\Dl$ such that $q \in \Psi_3(p)$, then
$p,q \in \Psi_2^2(r)$.\qed
\end{lem}

\begin{lem} \label{2}
Let $\{p,q\}$ be an edge of $\Dl$. Then exactly one of the following
holds:
\begin{itemize}
\item[{\rm (i)}] $|\Dl(p) \cap \Dl(q)|=(24t+4)(t^2+t+1)+1$ (which is
$2101$, $365$ and $85$ for $t=4,2$ and $1$) and $\{p,q\}$ is an edge of
$\Psi$;

\item[{\rm (ii)}] $|\Dl(p) \cap \Dl(q)|=8t^4+12t^3+4t^2+28t+5$
(which is $2997$, $391$ and $57$ for $t=4,2$ and $1$) and $\{p,q\}$ is
an edge of $\widehat \Psi$.\qed
\end{itemize}
\end{lem}

\vspace{0.5cm}

\vspace{1cm}

\unitlength 0.85mm
\linethickness{0.4pt}
\begin{picture}(169.29,142.03)
\put(82.00,19.37){\circle{8.00}}
\put(20.50,97.70){\oval(39.33,8.00)[]}
\put(82.50,134.03){\oval(39.33,8.00)[]}
\put(141.09,97.95){\oval(30.00,8.00)[]}
\put(82.00,23.37){\line(0,1){21.00}}
\put(82.00,19.37){\makebox(0,0)[cc]{\small 1}}
\put(21.17,97.70){\makebox(0,0)[cc]{\small $6(2t^2+1)(2t+1)$}}
\put(83.50,134.03){\makebox(0,0)[cc]{\small $96t^3(2t^2+1)(2t+1)$}}
\put(82.75,48.94){\makebox(0,0)[cc]{\small $16t(2t^2+1)(t^2+t+1)$}}
\put(141.09,98.28){\makebox(0,0)[cc]{\small $512t^6$}}
\put(81.58,48.61){\oval(45.00,8.00)[]}
\put(82.09,129.95){\line(0,-1){77.50}}
\put(62.50,52.45){\line(0,0){0.00}}
\put(37.92,93.70){\line(3,-5){24.50}}
\put(101.25,52.45){\line(2,3){27.50}}
\put(36.67,101.61){\line(1,1){28.75}}
\put(102.33,25.70){\makebox(0,0)[cc]{\scriptsize $16t(2t^2+1)(t^2+t+1)$}}
\put(86.33,42.03){\makebox(0,0)[cc]{\scriptsize 1}}
\put(34.33,48.70){\makebox(0,0)[cc]{\scriptsize $(10t-4)+4t^2(2t^2+3t+1)$}}
\put(97.34,77.37){\makebox(0,0)[cc]{\scriptsize $12t^2(2t^2+3t+1)$}}
\put(99.00,111.70){\makebox(0,0)[cc]{\scriptsize $2(t+1)(t^2+t+1)$}}
\put(113.67,60.70){\makebox(0,0)[cc]{\scriptsize $32t^5$}}
\put(140.00,82.37){\makebox(0,0)[cc]{\scriptsize $(2t^2+1)(t^2+t+1)$}}
\put(142.21,119.92){\makebox(0,0)[cc]{\scriptsize $3(2t^2+1)(2t+1)(t^2+t+1)$}}
\put(115.33,129.37){\makebox(0,0)[cc]{\scriptsize $16t^3(t^2+t+1)$}}
\put(48.00,128.37){\makebox(0,0)[cc]{\scriptsize $2(t^2+t+1)$}}
\put(32.66,111.70){\makebox(0,0)[cc]{\scriptsize $32t^3(t^2+t+1)$}}
\put(9.00,106.70){\makebox(0,0)[cc]{\scriptsize $8t(t^2+t+1)$}}
\put(29.00,86.37){\makebox(0,0)[cc]{\scriptsize $8t(t^2+t+1)$}}
\put(82.33,142.03){\makebox(0,0)[cc]{\scriptsize $2(8t^3+7t-2)(t^2+t+1)$}}
\put(48.34,58.37){\makebox(0,0)[cc]{\scriptsize $3(2t+1)$}}
\put(156.69,106.03){\makebox(0,0)[cc]{\scriptsize $-2t^2+6t-4$}}
\put(158.35,110.70){\makebox(0,0)[cc]{\scriptsize $20t^5+12t^4+22t^3-$}}
\put(94.67,19.37){\makebox(0,0)[cc]{\small $\{p\}$}}
\put(113.33,48.70){\makebox(0,0)[cc]{\small $\Psi_2^2(p)$}}
\put(51.00,97.70){\makebox(0,0)[cc]{\small $\Psi_1(p)$}}
\put(169.29,98.03){\makebox(0,0)[cc]{\small $\Psi_3(p)$}}
\put(110.21,139.21){\makebox(0,0)[cc]{\small $\Psi_2^4(p)$}}
\multiput(128.67,102.03)(-0.12,0.12){237}{\line(-1,0){0.12}}
\multiput(100.33,130.37)(0.12,-0.12){117}{\line(0,-1){0.12}}
\multiput(114.33,116.37)(-0.12,0.12){59}{\line(0,1){0.12}}
\put(81.85,10.00){\makebox(0,0)[cc]{Fig. 3}}
\end{picture}

The subgraph in $\Psi$ induced by (the points incident to) a symplecton
has the following suborbit diagram.

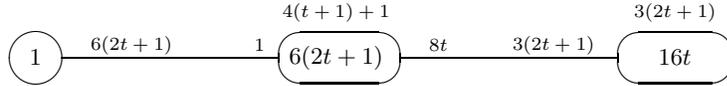
\unitlength 0.85mm
\linethickness{0.4pt}
\begin{picture}(134.00,16.00)
\put(25.00,9.00){\circle{8.00}}
\put(46.00,9.00){\line(1,0){17.00}}
\put(82.00,9.00){\line(1,0){17.00}}
\put(25.00,9.00){\makebox(0,0)[cc]{{\small 1}}}
\put(40.00,11.00){\makebox(0,0)[cc]{{\scriptsize $6(2t+1)$}}}
\put(60.00,11.00){\makebox(0,0)[cc]{{\scriptsize 1}}}
\put(88.00,11.00){\makebox(0,0)[cc]{{\scriptsize $8t$}}}
\put(106.00,11.00){\makebox(0,0)[cc]{{\scriptsize $3(2t+1)$}}}
\put(72.50,9.00){\oval(19.00,8.00)[]}
\put(72.00,9.00){\makebox(0,0)[cc]{\small $6(2t+1)$}}
\put(72.00,16.00){\makebox(0,0)[cc]{\scriptsize $4(t+1)+1$}}
\put(125.00,9.00){\oval(18.00,8.00)[]}
\put(125.00,9.00){\makebox(0,0)[cc]{\small $16t$}}
\put(125.00,16.00){\makebox(0,0)[cc]{\scriptsize $3(2t+1)$}}
\put(29.00,9.00){\line(1,0){17.00}}
\put(99.00,9.00){\line(1,0){17.00}}
\end{picture}

Comparing this diagram with that of $\Psi$ and $\widehat \Psi$ we obtain the
following

\begin{lem} \label{3}
Let $p \in \Dl$, $q \in \Psi_2^2(p)$ and let $\Upsilon$ be
the subgraph in $\Psi$ induced by (the points incident to)
the unique symplecton containing $p$ and $q$. Then
$$\Upsilon_1(p)=\Dl(q) \cap \Psi_1(p) ~~and~~$$
$$\Upsilon_2(p)=\{r \mid r \in \Psi_2^2(p), \Dl(r) \cap
\Psi_1(p)=\Upsilon_1(p)\}.$$
\end{lem}

\begin{lem} \label{3.5}
The automorphism group of $\Dl=$$^t\Dl$ is isomorphic to $^tF$ $($which is $^2E_6(2):S_3$, $F_4(2)$
or $\Omega_8^+(2):S_3$ for $t=4,2$ and $1$, respectively$)$.
\end{lem} 

\pf By (\ref{2}) the automorphism group of $\Dl$ does not mix the edges
of $\Psi$ with the edges of $\widehat \Psi$ and by (\ref{3}) the set of
symplecta of $\cF$ can be reconstructed from $\Dl$. The lines and
planes of $\cF$ can be reconstructed as the non-empty intersections of
two or more symplecta; the incidence relation in $\cF$ is via
inclusion. Hence every automorphism of $\Dl$ induces a uniquely
determined automorphism of $\cF$ and the result follows.\qed

\medskip

\section{$c$-extensions of $\cF$} \label{s3}
Thoughout this section $\cE=$$^t\cE$ is a geometry with
diagram
$c.F_4(t)$, where $t=4,2$ or $1$, such that the residue of an element
of type 1 is isomorphic to the building $^t\cF$ and the conditions $(a)$,
$(b)$ formulated in Section 1 are satisfied. Since these conditions
are equivalent to conditions  $(a^\prime)$ and $(b^\prime)$, we have the
following

\begin{lem} \label{5}
The graph $\Gm(\cE)$ is locally $\Dl$.\qed
\end{lem}

We will see that a converse of this lemma holds.

\medskip

Let $\Dl=$$^t\Dl$ for $t=4,2,1$ and let $\Gm$ be a graph which is locally
$\Dl$. For every $x \in \Gm$ let us fix a bijection
$$i_x: \Gm(x) \to \Dl$$
which induces an isomorphism of the subgraph in $\Gm$ induced by
$\Gm(x)$ onto $\Dl$ (so that if $y \in \Gm(x)$ then $i_x(y)$ is an
involution in $^tF$). The following lemma is a direct consequence of
(\ref{2}).

\begin{lem} \label{6}
Let $\{x_1,x_2,x_3\}$ be a triangle in $\Gm$. Then one of the following
holds:
\begin{itemize}
\item[{\rm (i)}] $i_{x_j}(x_k) \in \Psi_1(i_{x_j}(x_l))$ for all
$\{j,k,l\}=\{1,2,3\}$;
\item[{\rm (ii)}] $i_{x_j}(x_k) \in \Psi_2^2(i_{x_j}(x_l))$ for all
$\{j,k,l\}=\{1,2,3\}$.\qed
\end{itemize}
\end{lem}

The triangles in $\Gm$ as in (\ref{6} (i)) will be called {\it short}
while those in (\ref{6} (ii)) will be called {\it long}.

\begin{lem} \label{7}
[7]
Let $\Gm$ be a graph which is locally $\Dl$. Then $\Gm=\Gm(\cE)$,
for a $c.F_4(t)$-geometry $\cE$ as above.
\end{lem}

\pf Let $x$ be a vertex of $\Gm$ and $P$ be the point-set of a
symplecton in $\cF$. Then the set
$$Q_{x,P}=\{x\} \cup \{y \in \Gm(x) \mid i_x(y) \in P\}$$
is a complete subgraph of size $28t+8$ in $\Gm$. It is easy to deduce
using (\ref{3}) and (\ref{6}) that for every $y \in Q_{x,P}$ we have
$Q_{x,P}=Q_{y,R}$ for a symplecton $R$ in $\cF$.

Define $\cE$ to be a rank 5 geometry in which the elements of type 5
are the subgraphs $Q_{x,P}$ taken for all $x \in \Gm$ and all
symplecta $P$ in $\cF$; the remaining elements of $\cE$ are the
non-empty intersections of two or more such subgraphs;
the elements of type 1, 2, 3 and 4 are the intersections
of size 1, 2, 4 and 8, respectively. The incidence is via inclusion.
Then it is immediately seen that $\cE$ has the required properties and
that $\Gamma$ is the graph of $\cE$.\qed

\medskip

In what follows the elements of type 1, 2, 3, 4 and 5 in $\cE$ will be
identified with the corresponding complete subgraphs in $\Gm=\Gm(\cE)$
as in (\ref{7}) of size 1, 2, 4, 8 and $28t+8$, respectively.

\medskip

By (\ref{7}) it is clear that $\cE$ satisfies the
intersection property (cf. \cite{Pasi94}) and hence by \cite{CP92} the residue $\cH$ of an
element of type 5 in $\cE$  is a standard quotient of an affine polar
space of type $C_4$. On the other hand by condition $(b)$ in the considered situation
$\cH$ is a 1-point extension of its $C_3$-residue. This gives the following.

\begin{lem} \label{8}
Let $^t\cH$ be the residue of an element of type $5$ in $\cE$. Then $^t\cH$
is flag-transitive with the automorphism group isomorphic to
$$Sp_8(2),~~2^6:Sp_6(2),~~Sp_6(2)$$
for $t=4,2,1$, respectively, and acting doubly transitively on the set
of elements of type $1$ in $^t\cH$. Furthermore $^1\cH$ is a
subgeometry in $^2\cH$ and the latter is a subgeometry in $^4\cH$.\qed
\end{lem}

\begin{lem} \label{9.5}
Let $L=\{x,z,x_1,z_1\}$ be an element of type $3$ in $\cE$ and $y$
be an element of type $1$, not in $L$, but incident with $L$ to a
common element $Q$ of type $5$ and $N=L \cup \{y\}$. Then exactly one of the following
holds: 

\begin{itemize}
\item[{\rm (i)}] $N$ is contained in an element of type $4$ and every
triangle in $N$ is short;

\item[{\rm (ii)}] $N$ is not contained in an element of type $4$, there
are exactly two short triangles in $N$, say $T_1$ and $T_2$, which contain
$y$ and $N=T_1 \cup T_2$.

\end{itemize}
\end{lem}

\pf Consider $\cF={\rm res}_{\cE}(x)$. Then $M:=L
\setminus \{x\}$ is a line and $y$ is a point, both contained in the
symplecton $Q \setminus \{x\}$. Then by the basic property of polar spaces, either 
$y$ is collinear to every point in $M$, in which case $M \cup \{y\}$ is
contained in a plane of $\cF$ and we have (i); or $y$ is collinear to a
unique point, say $z$ in $M$, in which case $\{y,x,z\}$ is short and the remaining two
triangles in $N$ containing $\{y,x\}$ are long. Since a similar
assertion holds in ${\rm res}_{\cE}(x_1)$, the triangle $\{y,x_1,z_1\}$
must also be short.\qed

\medskip

>From (\ref{9.5}) we deduce that the short triangles contained in an
element of type 5 define a 2-graph, which is equivalent to the following 

\begin{cor} \label{9.75}
If $N$ is a $4$-element subset of $($the set of elements of type $1$
in$)$ an element of type $5$ in $\cF$, then the number of short
triangles in $N$ is even.\qed
\end{cor}
 
The following lemma is nothing but a more detailed version of (\ref{9.75}).

\begin{lem} \label{9.9}
Let $\{x,z\}$ and $\{u,v\}$ be disjoint elements of type $2$ contained
in an element of type $5$. Then exactly one of the following holds:

\begin{itemize}

\item[{\rm (i)}] $\{x,z,u,v\}$ is contained in an element of type $4$;

\item[{\rm (ii)}] both $\{x,z,u\}$ and $\{x,z,v\}$ are short and both
$\{x,u,v\}$ and $\{z,u,v\}$ are long;

\item[{\rm (iii)}] exactly one of $\{x,z,u\}$ and $\{x,z,v\}$ is short
and exactly one of $\{x,u,v\}$ and $\{z,u,v\}$ is short;

\item[{\rm (iv)}] both $\{x,z,u\}$ and $\{x,z,v\}$ are long and
$\{x,u,v\}$ and $\{z,u,v\}$ are of the same type (both short or both
long).\qed 

\end{itemize}
\end{lem}

The following result can alternatively be checked 
in the residue $^t\cH$, known by (\ref{8}).

\begin{lem} \label{9}
Let $X=\{x,z,u_1,u_2\}$ and $Y=\{x,z,v_1,v_2\}$ be distinct elements of
type $3$ in $\cE$ incident to the common element $\{x,z\}$ of type $2$ and
contained in a common element of type $5$. Then the symmetric
difference $Z=\{u_1,u_2,v_1,v_2\}$ of $X$ and $Y$ is also an element of
type $3$ in $\cE$.
\end{lem}

\pf If $X$ and $Y$ are contained in an element $N$ of type $4$, then
the result is immediate from the structure of ${\rm res}_{\cE}^-(N)$,
which is the $3$-dimensional $GF(2)$-space. Suppose that $X$ and $Y$
are not in an element of type $4$. Then applying (\ref{9.5}) to
$\{x,u_1,u_2,v_1\}$, we observe that the triangle $\{u_1,u_2,v_1\}$ is
short. By the obvious symmetry every triangle in $Z$ is short and hence
$Z$ is an element of type $3$.\qed

\section{The examples} \label{s4}
By (\ref{5}) and (\ref{7}) there is a natural bijection between the
$c.F_4(t)$-geometries satisfying $(a)$ and $(b)$ and the graphs which
are locally $^t\Dl$. Thus in order to construct a $c.F_4(t)$-geometry
it is sufficient to construct a graph which is locally $^t\Dl$.

\medskip

It is known (cf. Lemma~5.10.6 in \cite{I99}) that the Baby Monster
graph is locally $^4\Dl$ which gives the geometry $\cE(B)$ as in
\cite{B85}.

\medskip

For $t=4,2,1$ put $n(t)=|$$^t\Dl|$. The following well-known
result comes, for instance, from information
on the maximal 2-locals in $^tF$ (compare \cite{ATLAS}).

\begin{lem} \label{00}
Up to
conjugation in the automorphism group of $^tF$ every transitive
action of $^tF$ of degree $n(t)$ is similar to the action on the
vertex-set of $^t\Dl$.\qed
\end{lem}

Consider the action of $^2E_6(2):2$ on the cosets of $F_4(2) \times 2$. The
subdegrees of this action are calculated in \cite{L93}. In fact this
action is similar to the action of $O^2(B(a))$ on $\Th_3(a)$ where $B$ is
the Baby Monster and $\Th$ is the Baby Monster graph. This
observation enables one to deduce the subdegrees from the structure
constants of the centralizer algebra of $B$ acting on $\Th$ (cf. p. 128
in \cite{PS97}).

\begin{lem} \label{aa}
The action of $^2E_6(2):2$ on the cosets of $F_4(2) \times 2$ has rank
$4$ with subdegrees $1$, $n(2)=(2^4+1)(2^{12}-1)$,
$2^4(2^8-1)(2^{12}-1)$ and $2^{12}(2^4+1)(2^8+2^4+1)$.\qed
\end{lem}

Consider the orbital $\Gm$ of valency $n(2)$ of the action in
(\ref{aa}). Comparing the subdegrees in (\ref{aa}) with the subdegrees of the
action of $F_4(2)$ on $^2\Dl$ and since the action in (\ref{aa})
is by conjugation on a class of outer involutions in  $^2E_6(2):2$ it
is easy to see that $\Gm$ is locally $^2\Dl$ which gives the geometry
$\cE($$^2E_6(2))$.

\medskip

Similarly \cite{L93} gives the following.

\begin{lem} \label{bb}
The action of $E_6(2):2$ on $F_4(2) \times 2$ has rank $6$ with
subdegrees $1$, $n(2)=(2^4+1)(2^{12}-1)$,
$2^4(2^8-1)(2^{12}-1)$, $2^{12}(2^{12}-1)$ and
$2^{12}(2^4-1)(2^8-1)$ (twice).\qed
\end{lem}

Again we observe that the orbital of valency $n(2)$ is the graph of a
$c.F_4(2)$-geometry $\cE(E_6(2))$.

\medskip

The following result follows from Table 2 in \cite{CC88}

\begin{lem} \label{bc}
Let $V$ be a $26$-dimensional $GF(2)$-module of $F \cong F_4(2)$. Then
the orbit lengths of $F$ on the set of vectors of $V$ are the
following: $1$, $n(2)=(2^4+1)(2^{12}-1)$, $2^8(2^8+2^4+1)$,
$2^4(2^{12}-1)(2^8-1)$, $2^8(2^{12}-1)(2^4+1)$, $2^{12}(2^8-1)(2^4-1)$
and $2^{12}(2^{12}-1)$.\qed
\end{lem}

The orbit lengths in (\ref{bc}) are the subdegrees of the natural
action on $V$ of the semidirect product $V:F \cong 2^{26}:F_4(2)$ and
similarly to the above case we observe that the orbital of valency
$n(2)$ of this action is locally $^2\Dl$ and hence it is the graph of a
$c.F_4(2)$-geometry $\cE(2^{26}:F_4(2))$.

\medskip

>From Lemma 2.17.1 in \cite{ILLSS} or p. 112 in \cite{PS97}
we obtain the suborbits of the action of $Fi_{22}:2$ on the cosets of
$\Omega_8^+(2).S_3 \times 2$. Notice, that this action is similar to
the action of $B(a) \cap B(b)$ on $\Th_3(a) \cap \Th_3(b)$ for $b \in
\Th_2^3(a)$ (compare the intersection number $61\,776$ in the top matrix on p. 128 in
\cite{PS97}).

\begin{lem} \label{cc}
The action of $Fi_{22}:2$ on the cosets of $\Omega_8^+(2).S_3 \times 2$ has rank
$4$ with subdegrees $1$, $n(1)=1\,575$, $22\,400$, $37\,800$.\qed
\end{lem}

It is immediate from the structure constants of the action from
(\ref{cc}) given in \cite{ILLSS} that the orbital of valency $n(1)$ is
locally $^1\Dl$ and hence it is the graph of a $c.F_4(1)$-geometry
$\cE(Fi_{22})$.

\medskip

Let $^2H \cong $$^2E_6(2):2$, $^1H \cong Fi_{22}:2$ and $\Omega$ be the
set of cosets of $^2K \cong F_4(2) \times 2$ and $^1K \cong
\Omega_8^+(2).S_3 \times 2$, respectively. Then the elements of
$\Omega$ are identified with a conjugacy class of outer involutions in
$^tH$ so that the action is by conjugation. Let $d \in \Omega$ and
$\{d\}$, $\Omega_1(d)$, $\Omega_2^3(d)$, $\Omega_2^4(d)$ be the orbits
of $^tK$ on $\Omega$ (compare (\ref{aa}) and (\ref{cc})). Then under a
suitable numbering of the orbits we have $o(de)=2,3,4$ for $e \in
\Omega_1(d),\Omega_2^3(d),\Omega_2^4(d)$, respectively (so that
$^tK=C_{^tH}(d)$ and $|\Omega_1(d)|=n(t)$). For $t=2,1$ let
$^t \widetilde H$ be the unique non-split extension of $^tH$ by a normal
subgroup $Z$ of order 3 (recall that the 3-part of the Schur multiplier of $O^2(^tH)$
is of order 3). Let $\widetilde d$ be an involution in the preimage of
$\la d \ra$ in $^t\widetilde H$. Then $\widetilde d$ inverts $Z$
and hence $C_{^t\widetilde H}(\widetilde d)$ maps isomorphically onto
(and will be identified with) $^tK$. Let $\widetilde \Omega$ be 
the set of cosets of $^tK$ in $^t\widetilde H$. Notice that
$\widetilde \Omega$ is a conjugacy class of involutions in $^t
\widetilde H$ containing $\widetilde d$. It is easy to see that
the preimage of $\Omega_1(d)$ in $\widetilde \Omega$ splits into two
$^tK$-orbits with lengths $n(t)$ and $2 \cdot n(t)$. Let $\widetilde \Gm$ be
the orbital of valency $n(t)$. We claim that $\widetilde \Gm$ is
locally $^t\Dl$. Let $\widetilde e, \widetilde f \in \widetilde
\Gm(\widetilde d)$, $e,f$ be their respective images in $\Omega$ and
suppose that $o(ef)=2$. Since $\widetilde d$ is an involution which
commutes with $\widetilde e$ and $\widetilde f$, we have
$o(\widetilde e \widetilde f)=o((\widetilde e\widetilde d)(\widetilde
f\widetilde d))$. Since $Z$ commutes with the product of any two
involutions from $\widetilde \Omega$ the claim follows and we have the
following.

\begin{lem} \label{dd}
The geometries $\cE($$^2E_6(2))$ and $\cE(Fi_{22})$ possess $3$-fold
covers $\cE(3 \cdot $$^2E_6(2))$ and $\cE(3 \cdot Fi_{22})$, respectively.\qed
\end{lem}

The subdegrees of $3 \cdot $$^2E_6(2)$ on the cosets of
$F_4(2)$ are given in Table 3 in \cite{L93}.

\medskip

To the end of this section we discuss some $c.F_4(t)$-geometries which
do not satisfy condition $(b)$. One class of such examples can be constructed
as follows. Consider the action of the Baby Monster
$B$ on the cosets of a subgroup $2 \cdot $$^2E_6(2)$. Then one
of the orbitals of valency $n(4)$ with respect to the action (the graph
$\widetilde \Dl_2$ on p. 98 in \cite{I94}) is the
graph of a $c.F_4(4)$-geometry which is a double cover of $\cE(B)$. The
residue of an element of type 5 in this cover is a double cover of the
complete graph on 120 vertices. The preimages of $\cE(^2E_6(2))$- and
$\cE(Fi_{22})$-subgeometries in $\cE(B)$ are also proper double covers.
One can construct analogous double covers of $\cE(3 \cdot Fi_{22})$ and
$\cE(3 \cdot $$^2E_6(2))$.

Let $\cF=$$^t\cF$ be an $F_4$-building with 3 points per a line and let
$(R,\varphi)$ be a representation of $\cF$, which means that $R$ is a
group and $\varphi$ is the a mapping $\varphi:p \mapsto z_p$ of the
point-set of $\cF$ into $R$ such that (i) $R$ is generated by the image
of $\varphi$; (ii) $z_p^2=1$ for every point $p$; (iii) $z_pz_qz_r=1$
whenever $\{p,q,r\}$ is a line. For an arbitrary
element $u \in \cF$ define $\varphi(u)$ to be the subgroup of $R$
generated by the elements $z_p$ taken for all points $p$ incident to
$u$. Suppose that $(R,\varphi)$ is separable which means that 
the mapping $u \mapsto \varphi(u)$ of $\cF$ into the
set of subgroups of $R$ is injective. Then there is a standard
procedure which enables one to construct a $c$-extension $\cE$ of $\cF$:
the elements of $\cE$ are the elements of $R$ and all cosets of the
subgroups $\varphi(u)$ taken for all $u \in \cF$; the incidence is via
inclusion. Notice that $\cE(2^{26}:F_4(2))$ can be constructed along
these lines. 

\medskip

If $F$ is the automorphism group of $\cF$ then 
$(F^\infty,\psi)$ is a separable representation of $\cF$
where $\psi$ is the identity mapping. In addition for each $t=4,2$ and
$1$ the geometry $\cF$ possesses an abelian representation. By this reason
the classification of all (even flag-transitive) $c.F_4(t)$-geometries
seems to be far too difficult.

\section{A subgraph $\Xi$} \label{s5}
We follow the notation introduced in Section~\ref{s3} so that $\cE$ is
a $c.F_4(t)$-geometry for $t=4,2,1$ whose elements of type 1, 2, 3, 4
and 5 are identified with the corresponding complete subgraphs of size
1, 2, 4, 8 and $28t+8$ in the graph $\Gm=\Gm(\cE)$, which is locally
$^t\Dl$.

\medskip

Let $\widetilde \Xi=$$^t\widetilde \Xi$ be the graph on the set of elements of type 2 in $\cE$
in which two such elements are adjacent if they are incident to a
common element of type 3 but not to a common element of type 1
(equivalently if their union is an element of type 3). Let $e=\{x,y\}$
be an element of type 2 in $\cE$ (which is an edge of $\Gm$ and a
vertex of $\widetilde \Xi$). Then the set
$\widetilde \Xi(e)$ of neighbours of $e$ in $\widetilde \Xi$ is in the natural bijection
with the set of elements of type 3 incident to $e$ via the
correspondence
$$\psi:e_1 \to e \cup e_1$$
for $e_1 \in \widetilde \Xi(e)$. Let $\Lm$ be the graph on the set of elements of
type 3 incident to $e$ in which two of them are adjacent if incident to
a common element of type 4. Then $\Lm$ is the collinearity graph of the
residual dual polar space ${\rm res}^+_{\cE}(e)$ with the diagram
$\node_2\darc\node_t\arc\node_t$ and the suborbit diagram of $\Lm$ is

\medskip

\noindent
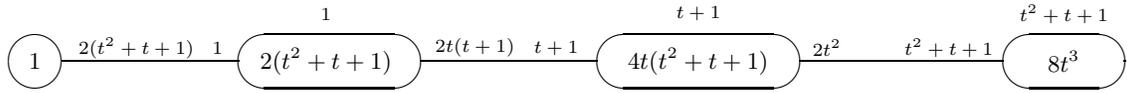
\unitlength 0.90mm
\linethickness{0.4pt}
\begin{picture}(166.00,12.00)
\put(5.00,5.00){\circle{8.00}}
\put(18.00,5.00){\line(1,0){17.00}}
\put(62.00,5.00){\line(1,0){17.00}}
\put(5.00,5.00){\makebox(0,0)[cc]{{\small 1}}}
\put(20.00,7.00){\makebox(0,0)[cc]{{\scriptsize $2(t^2+t+1)$}}}
\put(32.00,7.00){\makebox(0,0)[cc]{{\scriptsize 1}}}
\put(70.00,7.00){\makebox(0,0)[cc]{{\scriptsize $2t(t+1)$}}}
\put(82.00,7.00){\makebox(0,0)[cc]{{\scriptsize $t+1$}}}
\put(48.00,5.00){\makebox(0,0)[cc]{\small $2(t^2+t+1)$}}
\put(48.00,12.00){\makebox(0,0)[cc]{\scriptsize $1$}}
\put(103.00,5.00){\makebox(0,0)[cc]{\small $4t(t^2+t+1)$}}
\put(103.00,12.00){\makebox(0,0)[cc]{\scriptsize $t+1$}}
\put(9.00,5.00){\line(1,0){17.00}}
\put(71.00,5.00){\line(1,0){17.00}}
\put(118.00,5.00){\line(1,0){17.00}}
\put(122.00,7.00){\makebox(0,0)[cc]{{\scriptsize $2t^2$}}}
\put(140.00,7.00){\makebox(0,0)[cc]{{\scriptsize $t^2+t+1$}}}
\put(157.00,5.00){\oval(18.00,8.00)[]}
\put(157.00,5.00){\makebox(0,0)[cc]{\small $8t^3$}}
\put(157.00,12.00){\makebox(0,0)[cc]{\scriptsize $t^2+t+1$}}
\put(131.00,5.00){\line(1,0){17.00}}
\put(48.50,5.00){\oval(27.00,8.00)[]}
\put(103.00,5.00){\oval(30.00,8.00)[]}
\end{picture}

\bigskip

Let $e_1,e_2 \in \widetilde \Xi(e)$, where
$e_1=\{x_1,y_1\}$ and $e_2=\{x_2,y_2\}$. If $\psi(e_1)$ and $\psi(e_2)$ are at distance 1 or 2
in $\Lm$, then they are incident to a common element of type 5 in $\cE$
and by (\ref{9}) $e_1 \cup e_2$ is an element of type 3 which means
that $e_1$ and $e_2$ are adjacent in $\widetilde \Xi$. On the other hand, if
$\psi(e_1)$, $\psi(e_2)$ are at distance 3 in $\Lm$ then $i_{x}(x_2)
\in \Psi_2^4(i_x(x_1))$, in which case $e_1$ and $e_2$ are certainly
not adjacent in $\widetilde \Xi$. This implies the following.

\begin{lem} \label{10}
Locally $\widetilde \Xi$ is the distance $1$-or-$2$ graph of $\Lm$.\qed
\end{lem}

Now by Theorem 1.3 (vi), (vii), (ix) in \cite{C99} we obtain the
following.

\begin{prop} \label{11}
Let $e$ be an element of type $2$ in $\cE$ and $\widetilde \Xi^e$ be the connected component of
$\widetilde \Xi$ containing $e$. Then
\begin{itemize}
\item[{\rm (i)}] $^1\widetilde \Xi^e$ has $40$ vertices, its
automorphism group is isomorphic to
$U_4(2):2$ and the suborbit diagram is

\medskip

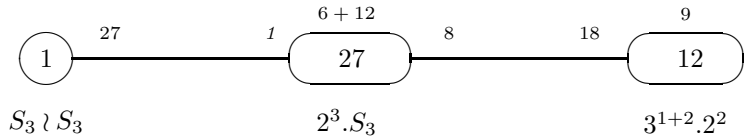
\unitlength 0.85mm
\linethickness{0.4pt}
\begin{picture}(134.00,19.00)
\put(25.00,12.00){\circle{8.00}}
\put(46.00,12.00){\line(1,0){17.00}}
\put(82.00,12.00){\line(1,0){17.00}}
\put(25.00,12.00){\makebox(0,0)[cc]{{$1$}}}
\put(35.00,16.00){\makebox(0,0)[cc]{{\scriptsize $27$}}}
\put(60.00,16.00){\makebox(0,0)[cc]{{\scriptsize 1}}}
\put(88.00,16.00){\makebox(0,0)[cc]{{\scriptsize $8$}}}
\put(110.00,16.00){\makebox(0,0)[cc]{{\scriptsize $18$}}}
\put(72.50,12.00){\oval(19.00,8.00)[]}
\put(72.00,12.00){\makebox(0,0)[cc]{ $27$}}
\put(72.00,19.00){\makebox(0,0)[cc]{\scriptsize $6+12$}}
\put(125.00,12.00){\oval(18.00,8.00)[]}
\put(125.00,12.00){\makebox(0,0)[cc]{ $12$}}
\put(125.00,19.00){\makebox(0,0)[cc]{\scriptsize $9$}}
\put(29.00,12.00){\line(1,0){17.00}}
\put(99.00,12.00){\line(1,0){17.00}}
\put(25.00,2.00){\makebox(0,0)[cc]{$S_3 \wr S_3$}}
\put(72.00,2.00){\makebox(0,0)[cc]{$2^3.S_3$}}
\put(125.00,2.00){\makebox(0,0)[cc]{$3^{1+2}.2^2$}}
\end{picture}

\item[{\rm (ii)}] $^2\widetilde \Xi^e$ has $256$ vertices, its
automorphism group is isomorphic
to $2^8:Sp_6(2)$ and the suborbit diagram is

\medskip

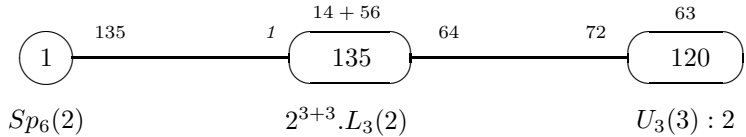
\unitlength 0.85mm
\linethickness{0.4pt}
\begin{picture}(134.00,20.00)
\put(25.00,13.00){\circle{8.00}}
\put(46.00,13.00){\line(1,0){17.00}}
\put(82.00,13.00){\line(1,0){17.00}}
\put(25.00,13.00){\makebox(0,0)[cc]{{$1$}}}
\put(35.00,17.00){\makebox(0,0)[cc]{{\scriptsize $135$}}}
\put(60.00,17.00){\makebox(0,0)[cc]{{\scriptsize 1}}}
\put(88.00,17.00){\makebox(0,0)[cc]{{\scriptsize $64$}}}
\put(111.00,17.00){\makebox(0,0)[cc]{{\scriptsize $72$}}}
\put(72.50,13.00){\oval(19.00,8.00)[]}
\put(72.00,13.00){\makebox(0,0)[cc]{ $135$}}
\put(72.00,20.00){\makebox(0,0)[cc]{\scriptsize $14+56$}}
\put(125.00,13.00){\oval(18.00,8.00)[]}
\put(125.00,13.00){\makebox(0,0)[cc]{ $120$}}
\put(125.00,20.00){\makebox(0,0)[cc]{\scriptsize $63$}}
\put(29.00,13.00){\line(1,0){17.00}}
\put(99.00,13.00){\line(1,0){17.00}}
\put(25.00,3.00){\makebox(0,0)[cc]{$Sp_6(2)$}}
\put(72.00,3.00){\makebox(0,0)[cc]{$2^{3+3}.L_3(2)$}}
\put(125.00,3.00){\makebox(0,0)[cc]{$U_3(3):2$}}
\end{picture}

\item[{\rm (iii)}] $^4\widetilde \Xi^e$ has $2300$ vertices, its
automorphism group is isomorphic
to $Co_2$ and the suborbit diagram is

\medskip

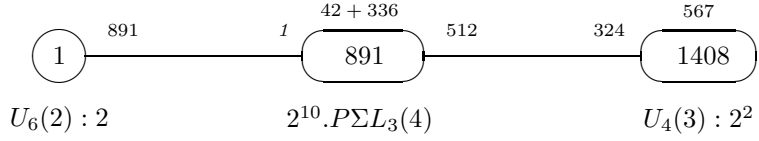
\unitlength 0.85mm
\linethickness{0.4pt}
\begin{picture}(134.00,19.00)
\put(25.00,12.00){\circle{8.00}}
\put(46.00,12.00){\line(1,0){17.00}}
\put(82.00,12.00){\line(1,0){17.00}}
\put(25.00,12.00){\makebox(0,0)[cc]{{$1$}}}
\put(35.00,16.00){\makebox(0,0)[cc]{{\scriptsize $891$}}}
\put(60.00,16.00){\makebox(0,0)[cc]{{\scriptsize 1}}}
\put(88.00,16.00){\makebox(0,0)[cc]{{\scriptsize $512$}}}
\put(111.00,16.00){\makebox(0,0)[cc]{{\scriptsize $324$}}}
\put(72.50,12.00){\oval(19.00,8.00)[]}
\put(72.00,12.00){\makebox(0,0)[cc]{ $891$}}
\put(72.00,19.00){\makebox(0,0)[cc]{\scriptsize $42+336$}}
\put(125.00,12.00){\oval(18.00,8.00)[]}
\put(125.00,12.00){\makebox(0,0)[cc]{ $1408$}}
\put(125.00,19.00){\makebox(0,0)[cc]{\scriptsize $567$}}
\put(29.00,12.00){\line(1,0){17.00}}
\put(99.00,12.00){\line(1,0){17.00}}
\put(25.00,2.00){\makebox(0,0)[cc]{$U_6(2):2$}}
\put(72.00,2.00){\makebox(0,0)[cc]{$2^{10}.P\Sigma L_3(4)$}}
\put(125.00,2.00){\makebox(0,0)[cc]{$U_4(3):2^2$}}
\end{picture}

\end{itemize}

\noindent
Furthermore, $^1\widetilde \Xi^e$ is a subgraph in $^2\widetilde \Xi^e$
and the latter is a
subgraph in $^4\widetilde \Xi^e$.\qed

\end{prop}

Let $\Xi^e$ be the subgraph in $\Gm$ induced by the elements of type 1
in $\cE$ incident
to the elements of type 2 in $\widetilde \Xi^e$. Since the diameter
of $\widetilde \Xi^e$ is 2, in view of the paragraph before
(\ref{10}), every vertex
in $\Xi^e$ is incident to exactly one element of type 2 in $\widetilde \Xi^e$
which gives the following.

\begin{lem} \label{12}
The subgraph $\Xi^e$ is the $2$-clique extension of $\widetilde \Xi^e$.\qed
\end{lem}

\medskip

\section{$\mu$-graphs of $D_8$-type} \label{s6}
In this and next sections for a pair $\{x,y\}$ of vertices at distance 2 in $\Gm$
we analyze the subgraph $\Gm(x,y)$ induced in $\Gm$ by the common
neighbours of $x$ and $y$ (the subgraphs $\Gm(x,y)$ are commonly known
as {\it $\mu$-graphs}).

\medskip

First we introduce some terminology. Recall that $\pi=(x,z,y)$ is
a 2-path in
$\Gm$ if $x,y \in \Gm(z)$ and $y \in \Gm_2(x)$. Such a $2$-path
$\pi$ is said to be of $D_{6}$- or $D_8$-type if
$$i_z(y) \in \Psi_3(i_z(x)) {\rm ~~or~~} i_z(y) \in \Psi_2^4(i_z(x)),$$
(equivalently if $\la i_z(x),i_z(y) \ra \cong D_6$ or $D_8$),
respectively. Clearly the paths $(x,z,y)$ and $(y,z,x)$ have the same
type.

\medskip

For $j=3$ or $4$ and a vertex $x$ of $\Gm$ let $\Gm_2^j(x)$
denote the set of vertices
$y$ at distance 2 from $x$ in $\Gm$ such that there is a path of
$D_{2j}$-type joining $x$ and $y$. Notice that $\Gm_2(x)=\Gm_2^3(x)
\cup \Gm_2^4(x)$ and {\it a priori}
the sets $\Gm_2^3(x)$ and $\Gm_2^4(x)$ might intersect.

\medskip

The next lemma follows from the fact that $\Gm$ is
locally $\Dl$ and from the suborbit diagrams of $\Psi$ and $\widehat \Psi$.

\begin{lem} \label{12.5}
Let $\pi=(x,z,y)$ be a $2$-path in $\Gm$. Then

\begin{itemize}

\item[{\rm (i)}] if $\pi$ is of $D_8$-type, then the valency of $z$ in
$\Gm(x,y)$ is $(2t+6)(t^2+t+1)+1$ (which is $295$, $71$ and $25$ for
$t=4,2$ and $1$);

\item[{\rm (ii)}] if $\pi$ is of $D_6$-type, then the valency of $z$ in
$\Gm(x,y)$ is $(2t^2+1)(t^2+t+1)$ (which is $693$, $63$ and $9$ for
$t=4,2$ and $1$).\qed

\end{itemize}

\end{lem}

\begin{lem} \label{13}
Let $\pi=(x,z,y)$ be a $2$-path of $D_8$-type. Then

\begin{itemize}

\item[{\rm (i)}] there is an element $e$ of type $2$ in $\cE$ incident
to $x$ such that $y,z \in \Xi^e$;

\item[{\rm (ii)}] every $2$-path contained in $\Xi^e$ is of $D_8$-type;

\item[{\rm (iii)}] the connected component of $\Gm(x,y)$ containing $z$
is contained in $\Xi^e$.

\end{itemize}

\end{lem}

\pf Since $\pi$ is of $D_8$-type by (\ref{1} (iv)) there is a (unique)
vertex $u$ adjacent to $z$, such that $i_z(x), i_z(y) \in
\Psi_1(i_z(u))$.  Then clearly $\pi \in \Xi^e$ where $e=\{z,u\}$ and
(i) follows. Now (ii) follows from the paragraph before (\ref{10}) and
the fact that the diameter of $\Xi^e$ is 2. In order to prove (iii)
it is sufficient to show that the valency
$k_1$ of $z$ in $\Gm(x,y) \cap \Xi^e$ is equal to the valency of $z$ in
$\Gm(x,y)$ given in (\ref{12.5} (i)).
By (\ref{10}) the graph $\widetilde \Xi$ is locally the distance
$1$-or-$2$ graph of $\Lm$. The
latter is a near hexagon with classical quads. This shows that for $f
\in \widetilde \Xi$, $g \in \widetilde \Xi_2(f)$ the subgraph in $\widetilde \Xi$ induced
by $\widetilde \Xi(f) \cap \widetilde \Xi(g)$ has valency
$k_2=(t+3)(t^2+t+1)$. Since $\Xi^e$ is the $2$-clique extension of
$\widetilde \Xi^e$ we have $k_1=2\times k_2+1$, hence (iii) follows.\qed

\medskip

\begin{lem} \label{13.5}
Let $\pi=(x,z,y)$ and $\sigma=(x,u,y)$ be $2$-paths in $\Gm$ such that
$z$ and $u$ are in the same connected component of $\Gm(x,y)$ and
$\delta=(z,v,w)$ be a $2$-path contained in $\Gm(x,y)$. Then $\pi$, $\sigma$ and
$\delta$ have the same type.
\end{lem}

\pf If one of the paths $\pi$, $\sigma$, $\delta$ is of $D_8$-type then
by (\ref{13}) all three paths are of $D_8$-type and hence the result.\qed

\medskip

Let $y \in \Gm_2^4(x)$, $e=\{x,z\}$ be an element of type 2 in $\cE$
such that $y \in \Xi^e$ (compare (\ref{13} (i)), $p=i_x(z )$, 
$\Omega=\{i_x(v) \mid v \in \Gm(x,y)\cap \Xi^e\}$ and $I$
be the setwise stabilizer of $\Omega$ in $F=$$^tF$. It is clear that $I
\le F(p)$.

\begin{lem} \label{14}
\begin{itemize}

\item[{\rm (i)}] $I$ contains $O_2(F(p))$;

\item[{\rm (ii)}] for every element $\alpha$ in $O^2(F(p))$ there is an
automorphism $\beta$ of $\Xi^e$ stabilizing $e$ such that the action of
$\alpha$ on the set of elements of type $3$ in $\cE$ incident to $e$
is similar to the action of $\beta$ on $\widetilde \Xi(e)$;

\item[{\rm (iii)}] $I \cong 2^{1+20}_+:U_4(3):2^2$, $2^{1+6+8}:U_3(3):2$ and
$2^{1+8}_+.3^{1+2}.2^2$ for $t=4,2$ and $1$, respectively.

\end{itemize}
\end{lem}

\pf An element $e_1$ of type 2 in $\cE$ is an orbit of $O_2(F(p))$ on $\Psi_1(p)$ if and only if
$e_1$ is the image under $i_x$ of an element from $\widetilde \Xi(e)$,
we obtain (i). Now (ii) and (iii) are immediate from (\ref{10}).\qed

\medskip

With $\Omega$ and $I$ as above put $\mu_1=|\Omega|$. Our next goal is to
describe the orbits of $I$ on $\Dl$.

\begin{lem} \label{15}
Let $L$ and $S$ be the set of lines and the set of symplecta in
$\cF$ incident to $p$, so that $|L|=3(2t^2+1)(2t+1)$ and
$|S|=(2t^2+1)(t^2+t+1)$. Then
\begin{itemize}
\item[{\rm (i)}] $I$ has two orbits, $L_1$ and $L_2$ on $L$ with
lengths $\mu_1/2$ and $|L|-\mu_1/2$, respectively;
\item[{\rm (ii)}] if $t=1$ or $2$ then $I$ acts transitively on $S$;
\item[{\rm (iii)}] if $t=4$ then $I$ has two orbits $S_1$ and $S_2$ on
$S$ with lengths $126$ and $567$, furthermore,
\begin{itemize}
\item[{\rm (a)}] the symplecta incident to a line from $L_1$ are
contained in $S_2$;
\item[{\rm (b)}] a line from $L_2$ is incident to $6$ symplecta from
$S_1$ and to $15$ symplecta from $S_2$.
\end{itemize}
\end{itemize}
\end{lem}

\pf The result follows from well-known properties of
$F(p)/O_2(F(p))$. For the case $t=4$ see Lemma 5.10.15 in \cite{I99}.\qed

\medskip

The next result is a direct consequence of (\ref{1} (vii)), (\ref{14} (i))
and (\ref{15}).

\begin{lem} \label{16}
The following assertions hold:
\begin{itemize}
\item[{\rm (i)}] $I$ has two orbits $\Omega_1=\Omega$ and
$\Omega_2$ on $\Psi_1(p)$;
\item[{\rm (ii)}] if $t=1,2$ then $I$ acts transitively on
$\Psi_2^2(p)$ while if $t=4$, $I$ has two orbits $\Omega_3$ and
$\Omega_4$ on $\Psi_2^2(p)$ with lengths $2^6 \cdot 126$ and $2^6 \cdot
567$, respectively;
\item[{\rm (iii)}] $I$ has two orbits $\Omega_5$ and $\Omega_6$ on
$\Psi_2^4(p)$ with lengths
$16t^3 \cdot |\Omega_1|$ and $16t^3 \cdot |\Omega_2|$;
\item[{\rm (iv)}] $I$ acts transitively on $\Omega_7=\Psi_3(p)$.\qed
\end{itemize}
\end{lem}

Notice that in terms of (\ref{15})
the orbits $\Omega_1$ and $\Omega_5$ correspond to $L_1$,
$\Omega_2$ and $\Omega_6$ correspond to $L_2$, $\Omega_3$
corresponds to $S_1$ and $\Omega_4$ corresponds to $S_2$.

\begin{lem} \label{17}
Let $u \in \Psi_1(p) \cup \Psi_2^2(p) \cup \Psi_2^4(p)$ for $t=1,2$ and
$u \in \Psi_1(p) \cup \Omega_4 \cup \Psi_2^4(p)$ for $t=4$. Then $u$ is adjacent
in $\Dl$ to a vertex from $\Omega$.
\end{lem}

\pf If $u \in \Psi_1(p)$, then the claim follows from the fact that the
subgraph in $\Dl$ induced by $\Psi_1(p)$ is
connected. A vertex from $\Omega$ is clearly adjacent in $\Dl$ to
vertices in $\Psi_2^2(p)$. Thus the claim follows for $t=1,2$ and $u \in
\Psi_2^2(p)$ as well as for $t=4$ and $u \in \Omega_4$ (compare (\ref{15}
(a)). Let $u \in \Psi_2^4(p)$ and $v$ be the unique vertex in $\Psi_1(p)
\cap \Psi_1(u)$ (compare the diagram of $\Psi$). Let $A$ and $B$ be the
orbits of $O_2(F(p))$ containing $u$ and $v$, respectively (so that
by (\ref{1} (vi)) we have $|A|=2^5 \cdot t^3$ and $
|B|=2$).
Then by (\ref{1} (vii)) the
stabilizer of $A$ in $F(p)$ coincides with the stabilizer of $B$ in
$F(p)$ and it is the stabilizer of the unique line $l$ of $\cF$
containing $p$ and $v$. Let $\cC$ be the set of orbits of $O_2(F(p))$
on $\Psi_1(p)$ containing vertices adjacent to $u$ in $\Dl$.
Since $|\Dl(p) \cap
\Dl(u)|=1+2(t^2+t+1)$ while the suborbits of the action of $F(p)$ on the set of
lines incident to $p$ are $1$, $2(t^2+t+1)$, $4t(t^2+t+1)$, $8t^3$
(compare the diagram of $\Lm$ before (\ref{10})), we conclude that
$|\cC|=1+2(t^2+t+1)$ and $\cC\setminus B$
correspond to the lines incident with $l$ to a common plane. Since
every vertex from $\Psi_1(p)$ is adjacent in $\Psi$ to a vertex from
$\Omega$, this completes the proof.\qed

\medskip

As a consequence of the above proof in view of (\ref{1} (v)) and (\ref{15}
(iii) (a)) we obtain the following.

\begin{cor} \label{cor}
If $u \in \Psi_3(p)$ for $t=1,2$ and $u \in \Omega_3 \cup \Psi_3(p)$
for $t=4$, then no vertex in $\Omega$ is adjacent to $u$ in $\Dl$.\qed
\end{cor}

\begin{lem} \label{18}
In terms introduced before $(\ref{14})$ if $w \in \Gm(x,y) \setminus 
\Xi^e$, then $(x,w,y)$ is of $D_6$-type and either $t=4$ or $i_x(w)
\in \Psi_3(p)$. 
\end{lem}

\pf By (\ref{13} (iii)) $\Gm(x,y) \cap \Xi^e$ is a connected component of
$\Gm(x,y)$ and hence by (\ref{17}) $i_x(w) \in \Psi_3(p)$ for $t=1,2$
and $i_x(w) \in \Psi_3(p) \cup \Omega_3$ for $t=4$. Suppose that
$(x,w,y)$ is of $D_8$-type. By
(\ref{13} (i), (iii)) there is an element $f=\{x,v\}$ of type 2 in
$\cE$ such that $y,w \in \Xi^f$ and $\Gm(x,y) \cap \Xi^f$ is
another connected component of $\Gm(x,y)$. Put $q=i_x(v)$ and
$\Omega'=\{i_x(s) \mid s \in \Gm(x,y) \cap
\Xi^f\}$. Recall that $\Omega'$ is the
union of some $\mu_1/2$ lines containing $q$ with $q$ removed. It is
immediate from the diagram of $\Psi$ that a line intersecting
$\Psi_3(p)$ contains two points from $\Psi_3(p)$ and one point from
$\Psi_2^4(p)$. Hence, assuming that $\Omega' \cap \Psi_3(p) \ne
\emptyset$, we must have $q \in \Psi_2^4(p)$. In this case let $r$ be
the unique vertex in $\Psi_1(p) \cap \Psi_1(q)$. Then by (\ref{17}) $r$
is adjacent in $\Dl$ to a vertex from $\Omega'$ and by (\ref{1.5})
$r$ is not adjacent to any vertex from $\Psi_3(p)$. Hence $\Omega' \cap
\Omega_3 \ne \emptyset$, in particular $t=4$. Let $\{q,s_1,s_2\}$ be a line
such that $s_1,s_2 \in \Omega'$ and $s_1 \in \Omega_3$. Then one can
see from Fig. 2 that $s_2 \not \in \Psi_3(p)$ and hence $s_2 \in \Omega_3$.
Comparing the diagram of $\Psi$ and the diagram of
the subgraph in $\Psi$ induced by a symplecton given before (\ref{3})
we observe that $s_1$ and $s_2$ are contained in the same symplecton,
containing $p$ and $q \in \Psi_1(p)$, 
in particular $\Omega' \cap \Psi_3(p)=\emptyset$.
By (\ref{15} (iii) (b)) $q$ is contained in exactly 6
symplecta intersecting $\Omega_3$ and each of these symplecta
contains just 16 lines containing $q$ and intersecting $\Psi_2^2(p)$. Since $6 \times 16 <
324$, we have a contradiction, hence $(x,w,y)$ must be of $D_6$-type.\qed

\medskip

The results in (\ref{11}), (\ref{13}) and (\ref{18}) can be summarized
as follows.

\begin{prop} \label{19}
Let $y \in \Gm_2^4(x)$. Then
\begin{itemize}
\item[{\rm (i)}] there is a unique element $e$ of type $2$ in $\cE$
incident to $x$ such that $y \in \Xi^e$;
\item[{\rm (ii)}] for any two elements $e$ and $f$ of type $2$ in $\cE$
the intersection $\Xi^e \cap \Xi^f$ induces in $\Gm$ a complete
subgraph;

\item[{\rm (iii)}] $|\Gm_2^4(x)|=|\Dl|\cdot |\Psi_2^4(p)|/\mu_1$
where $\mu_1=648$, $144$ and $36$ for $t=4,2$ and $1$.\qed
\end{itemize}
\end{prop}

\section{$\mu$-graphs of $D_6$-type} \label{s7}
In this section we analyze the common neighbours of $x$ and $y \in \Gm_2^3(x)$. Let
$(x,z,y)$ be a $2$-path of $D_6$-type, $p=i_z(x)$, $q=i_z(y)$ and
$^t\Upsilon=\Upsilon=\{i_z(u) \mid u \in \Gm(x,y) \cap \Gm(z)\}$. Then $q \in
\Psi_3(p)$ and by (\ref{1} (v)) $H=$$^tH:=F(p) \cap F(q)$ is a complement to
$O_2(F(p))$ in $F(p)$ isomorphic to $U_6(2).S_3$, $Sp_6(2)$ and $S_3
\wr S_3$ for $t=4,2$ and $1$.

\begin{lem} \label{20}
The following assertions hold:
\begin{itemize}
\item[{\rm (i)}] $|\Upsilon|=(2t^2+1)(t^2+t+1)$;
\item[{\rm (ii)}] $\Upsilon$ is the unique class of central involutions
in $H$ of size $|\Upsilon|$;
\item[{\rm (iii)}] $\Upsilon$ is a class of $3$-transpositions in $H$,
so that
\begin{itemize}

\item[{\rm (a)}] $^1\Upsilon$ is the complete $3$-partite graph
$K_{3 \times 3}$ with the suborbit diagram

\medskip

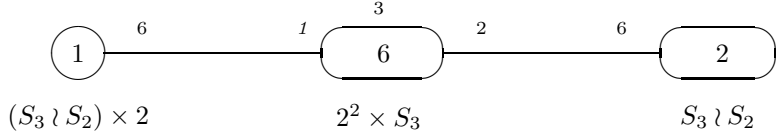
\unitlength 0.85mm
\linethickness{0.4pt}
\begin{picture}(134.00,19.00)
\put(25.00,12.00){\circle{8.00}}
\put(46.00,12.00){\line(1,0){17.00}}
\put(82.00,12.00){\line(1,0){17.00}}
\put(25.00,12.00){\makebox(0,0)[cc]{{$1$}}}
\put(35.00,16.00){\makebox(0,0)[cc]{{\scriptsize $6$}}}
\put(60.00,16.00){\makebox(0,0)[cc]{{\scriptsize 1}}}
\put(88.00,16.00){\makebox(0,0)[cc]{{\scriptsize $2$}}}
\put(110.00,16.00){\makebox(0,0)[cc]{{\scriptsize $6$}}}
\put(72.50,12.00){\oval(19.00,8.00)[]}
\put(72.00,12.00){\makebox(0,0)[cc]{ $6$}}
\put(72.00,19.00){\makebox(0,0)[cc]{\scriptsize $3$}}
\put(125.00,12.00){\oval(18.00,8.00)[]}
\put(125.00,12.00){\makebox(0,0)[cc]{ $2$}}
\put(29.00,12.00){\line(1,0){17.00}}
\put(99.00,12.00){\line(1,0){17.00}}
\put(25.00,2.00){\makebox(0,0)[cc]{$(S_3 \wr S_2)\times 2$}}
\put(72.00,2.00){\makebox(0,0)[cc]{$2^2 \times S_3$}}
\put(125.00,2.00){\makebox(0,0)[cc]{$S_3 \wr S_2$}}
\end{picture}

\item[{\rm (b)}] $^2\Upsilon$ is the collinearity graph of the polar
space $\cP(Sp_6(2))$ of $^2H \cong Sp_6(2)$ with the diagram
$\node_2\arc\node_2\darc\node_2$ and the suborbit diagram of
$^2\Upsilon$ is

\medskip

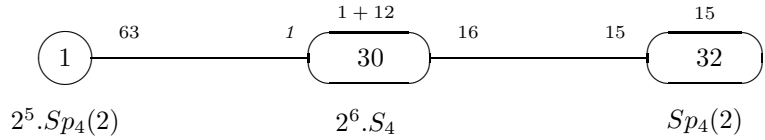
\unitlength 0.85mm
\linethickness{0.4pt}
\begin{picture}(134.00,20.00)
\put(25.00,13.00){\circle{8.00}}
\put(46.00,13.00){\line(1,0){17.00}}
\put(82.00,13.00){\line(1,0){17.00}}
\put(25.00,13.00){\makebox(0,0)[cc]{{$1$}}}
\put(35.00,17.00){\makebox(0,0)[cc]{{\scriptsize $63$}}}
\put(60.00,17.00){\makebox(0,0)[cc]{{\scriptsize 1}}}
\put(88.00,17.00){\makebox(0,0)[cc]{{\scriptsize $16$}}}
\put(111.00,17.00){\makebox(0,0)[cc]{{\scriptsize $15$}}}
\put(72.50,13.00){\oval(19.00,8.00)[]}
\put(72.00,13.00){\makebox(0,0)[cc]{ $30$}}
\put(72.00,20.00){\makebox(0,0)[cc]{\scriptsize $1+12$}}
\put(125.00,13.00){\oval(18.00,8.00)[]}
\put(125.00,13.00){\makebox(0,0)[cc]{ $32$}}
\put(125.00,20.00){\makebox(0,0)[cc]{\scriptsize $15$}}
\put(29.00,13.00){\line(1,0){17.00}}
\put(99.00,13.00){\line(1,0){17.00}}
\put(25.00,3.00){\makebox(0,0)[cc]{$2^5.Sp_4(2)$}}
\put(72.00,3.00){\makebox(0,0)[cc]{$2^{6}.S_4$}}
\put(125.00,3.00){\makebox(0,0)[cc]{$Sp_4(2)$}}
\end{picture}

\item[{\rm (c)}] $^4\Upsilon$ is the collinearity graph of the polar
space $\cP(U_6(2))$ of $^4H \cong U_6(2).S_3$ with the diagram
$\node_4\arc\node_4\darc\node_2$ and the suborbit diagram of $^4\Upsilon$ is

\medskip

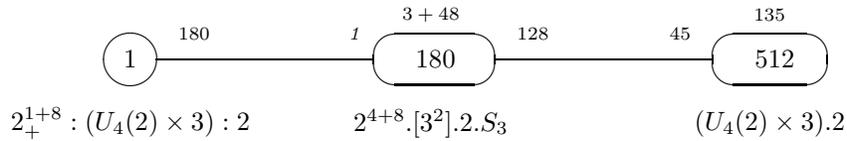
\unitlength 0.85mm
\linethickness{0.4pt}
\begin{picture}(134.00,19.00)
\put(25.00,12.00){\circle{8.00}}
\put(46.00,12.00){\line(1,0){17.00}}
\put(82.00,12.00){\line(1,0){17.00}}
\put(25.00,12.00){\makebox(0,0)[cc]{{$1$}}}
\put(35.00,16.00){\makebox(0,0)[cc]{{\scriptsize $180$}}}
\put(60.00,16.00){\makebox(0,0)[cc]{{\scriptsize 1}}}
\put(88.00,16.00){\makebox(0,0)[cc]{{\scriptsize $128$}}}
\put(111.00,16.00){\makebox(0,0)[cc]{{\scriptsize $45$}}}
\put(72.50,12.00){\oval(19.00,8.00)[]}
\put(72.00,12.00){\makebox(0,0)[cc]{ $180$}}
\put(72.00,19.00){\makebox(0,0)[cc]{\scriptsize $3+48$}}
\put(125.00,12.00){\oval(18.00,8.00)[]}
\put(125.00,12.00){\makebox(0,0)[cc]{ $512$}}
\put(125.00,19.00){\makebox(0,0)[cc]{\scriptsize $135$}}
\put(29.00,12.00){\line(1,0){17.00}}
\put(99.00,12.00){\line(1,0){17.00}}
\put(25.00,2.00){\makebox(0,0)[cc]{$2^{1+8}_+:(U_4(2)\times 3):2$}}
\put(72.00,2.00){\makebox(0,0)[cc]{$2^{4+8}.[3^2].2.S_3$}}
\put(125.00,2.00){\makebox(0,0)[cc]{$(U_4(2) \times 3).2$}}
\end{picture}

\end{itemize}

\end{itemize}

\end{lem}

\pf The assertion (i) is just (\ref{12.5} (ii)). By the definition,
$\Upsilon$ consists of all the involutions in $\Dl$ which commute with
$p$ and $q$. Hence $\Upsilon$ is contained in $H$ and it is closed
under conjugation by the elements of $H$. Hence $\Upsilon$ is a union
of some conjugacy classes of involutions in $H$ and since $|\Upsilon|$ is
odd, at least one of these classes consists of central involutions.
Now (ii) and (iii) follows from the well-known properties of $H$ (cf.
\cite{ATLAS}).\qed

\medskip

In terms introduced before (\ref{20}) let $\Pi$ be the connected
component of $\Gm(x,y)$ containing $z$ and $\Sg=\{i_x(u) \mid u \in
\Pi\}$.

\begin{lem} \label{21}
Let $u,v \in \Sg$. Then
\begin{itemize}
\item[{\rm (i)}] $v \in \Psi_2^2(u)$ if $u$
and $v$ are adjacent in $\Sg$ and $v \in \Psi_3(u)$ if $u$ and $v$ are
at distance $2$ in $\Sg$;
\item[{\rm (ii)}] the graph $\Sg$ is locally $\Upsilon$.
\end{itemize}
\end{lem}

\pf The assertion (ii) follows directly from the definition of $\Sg$
and $\Upsilon$ while (i) is immediate from (\ref{1.5}) and
(\ref{13.5}).\qed

\begin{lem} \label{22}
\begin{itemize}
\item[{\rm (i)}] $^1\Sg$ is the complete $4$-partite graph
$K_{4 \times 3}$ with the automorphism group $S_3 \wr S_4$ and the
suborbit diagram

\medskip

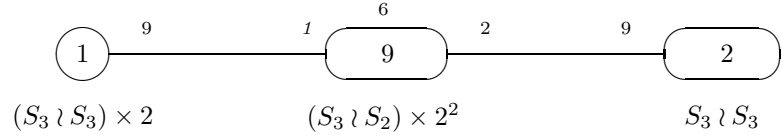
\unitlength 0.85mm
\linethickness{0.4pt}
\begin{picture}(148.00,19.00)
\put(39.00,12.00){\circle{8.00}}
\put(60.00,12.00){\line(1,0){17.00}}
\put(96.00,12.00){\line(1,0){17.00}}
\put(39.00,12.00){\makebox(0,0)[cc]{{$1$}}}
\put(49.00,16.00){\makebox(0,0)[cc]{{\scriptsize $9$}}}
\put(74.00,16.00){\makebox(0,0)[cc]{{\scriptsize 1}}}
\put(102.00,16.00){\makebox(0,0)[cc]{{\scriptsize $2$}}}
\put(124.00,16.00){\makebox(0,0)[cc]{{\scriptsize $9$}}}
\put(86.50,12.00){\oval(19.00,8.00)[]}
\put(86.00,12.00){\makebox(0,0)[cc]{ $9$}}
\put(86.00,19.00){\makebox(0,0)[cc]{\scriptsize $6$}}
\put(139.00,12.00){\oval(18.00,8.00)[]}
\put(139.00,12.00){\makebox(0,0)[cc]{ $2$}}
\put(43.00,12.00){\line(1,0){17.00}}
\put(113.00,12.00){\line(1,0){17.00}}
\put(39.00,2.00){\makebox(0,0)[cc]{$(S_3 \wr S_3) \times 2$}}
\put(86.00,2.00){\makebox(0,0)[cc]{$(S_3 \wr S_2) \times 2^2$}}
\put(139.00,2.00){\makebox(0,0)[cc]{$S_3 \wr S_3$}}
\end{picture}

\medskip

\item[{\rm (ii)}] $^2\Sg$ is isomorphic to one of the following
three graphs:

\begin{itemize}
\item[{\rm (a)}] the commuting graph $^2\Sg^a$ of
$3$-transpositions in $\Omega_8^-(2):2$
with the suborbit diagram

\medskip

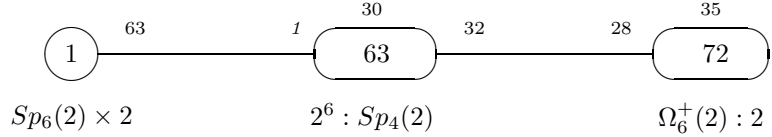
\unitlength 0.85mm
\linethickness{0.4pt}
\begin{picture}(142.00,20.00)
\put(33.00,13.00){\circle{8.00}}
\put(54.00,13.00){\line(1,0){17.00}}
\put(90.00,13.00){\line(1,0){17.00}}
\put(33.00,13.00){\makebox(0,0)[cc]{{$1$}}}
\put(43.00,17.00){\makebox(0,0)[cc]{{\scriptsize $63$}}}
\put(68.00,17.00){\makebox(0,0)[cc]{{\scriptsize 1}}}
\put(96.00,17.00){\makebox(0,0)[cc]{{\scriptsize $32$}}}
\put(119.00,17.00){\makebox(0,0)[cc]{{\scriptsize $28$}}}
\put(80.50,13.00){\oval(19.00,8.00)[]}
\put(80.00,13.00){\makebox(0,0)[cc]{ $63$}}
\put(80.00,20.00){\makebox(0,0)[cc]{\scriptsize $30$}}
\put(133.00,13.00){\oval(18.00,8.00)[]}
\put(133.00,13.00){\makebox(0,0)[cc]{ $72$}}
\put(133.00,20.00){\makebox(0,0)[cc]{\scriptsize $35$}}
\put(37.00,13.00){\line(1,0){17.00}}
\put(107.00,13.00){\line(1,0){17.00}}
\put(33.00,3.00){\makebox(0,0)[cc]{$Sp_6(2)\times 2$}}
\put(80.00,3.00){\makebox(0,0)[cc]{$2^6:Sp_4(2)$}}
\put(133.00,3.00){\makebox(0,0)[cc]{$\Omega_6^+(2):2$}}
\end{picture}

\medskip

\item[{\rm (b)}] the commuting graph $^2\Sg^b$ of
$3$-transpositions in $\Omega_8^+(2):2$ with the suborbit diagram

\medskip

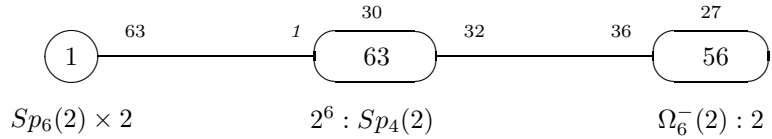
\unitlength 0.85mm
\linethickness{0.4pt}
\begin{picture}(142.00,20.00)
\put(33.00,13.00){\circle{8.00}}
\put(54.00,13.00){\line(1,0){17.00}}
\put(90.00,13.00){\line(1,0){17.00}}
\put(33.00,13.00){\makebox(0,0)[cc]{{$1$}}}
\put(43.00,17.00){\makebox(0,0)[cc]{{\scriptsize $63$}}}
\put(68.00,17.00){\makebox(0,0)[cc]{{\scriptsize 1}}}
\put(96.00,17.00){\makebox(0,0)[cc]{{\scriptsize $32$}}}
\put(119.00,17.00){\makebox(0,0)[cc]{{\scriptsize $36$}}}
\put(80.50,13.00){\oval(19.00,8.00)[]}
\put(80.00,13.00){\makebox(0,0)[cc]{ $63$}}
\put(80.00,20.00){\makebox(0,0)[cc]{\scriptsize $30$}}
\put(133.00,13.00){\oval(18.00,8.00)[]}
\put(133.00,13.00){\makebox(0,0)[cc]{ $56$}}
\put(133.00,20.00){\makebox(0,0)[cc]{\scriptsize $27$}}
\put(37.00,13.00){\line(1,0){17.00}}
\put(107.00,13.00){\line(1,0){17.00}}
\put(33.00,3.00){\makebox(0,0)[cc]{$Sp_6(2)\times 2$}}
\put(80.00,3.00){\makebox(0,0)[cc]{$2^6:Sp_4(2)$}}
\put(133.00,3.00){\makebox(0,0)[cc]{$\Omega_6^-(2):2$}}
\end{picture}

\medskip

\item[{\rm (c)}] the $2$-fold antipodal cover $^2\Sg^c$ of the complete graph
with the automorphism group $2^7:Sp_6(2)$ and the suborbit diagram

\medskip

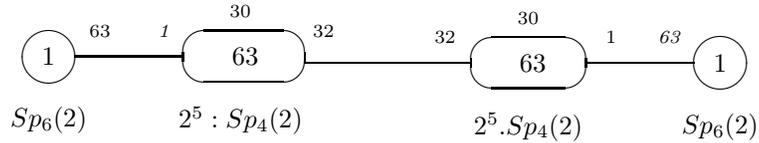
\unitlength 0.85mm
\linethickness{0.4pt}
\begin{picture}(145.00,21.00)
\put(36.00,14.00){\circle{8.00}}
\put(85.00,13.00){\line(1,0){17.00}}
\put(36.00,14.00){\makebox(0,0)[cc]{{$1$}}}
\put(44.00,18.00){\makebox(0,0)[cc]{{\scriptsize $63$}}}
\put(54.00,18.00){\makebox(0,0)[cc]{{\scriptsize 1}}}
\put(79.00,18.00){\makebox(0,0)[cc]{{\scriptsize $32$}}}
\put(98.00,17.00){\makebox(0,0)[cc]{{\scriptsize $32$}}}
\put(66.50,14.00){\oval(19.00,8.00)[]}
\put(66.00,14.00){\makebox(0,0)[cc]{ $63$}}
\put(66.00,21.00){\makebox(0,0)[cc]{\scriptsize $30$}}
\put(111.00,13.00){\oval(18.00,8.00)[]}
\put(111.00,13.00){\makebox(0,0)[cc]{ $63$}}
\put(111.00,20.00){\makebox(0,0)[cc]{\scriptsize $30$}}
\put(40.00,14.00){\line(1,0){17.00}}
\put(36.00,4.00){\makebox(0,0)[cc]{$Sp_6(2)$}}
\put(66.00,4.00){\makebox(0,0)[cc]{$2^5:Sp_4(2)$}}
\put(111.00,3.00){\makebox(0,0)[cc]{$2^5.Sp_4(2)$}}
\put(124.00,17.00){\makebox(0,0)[cc]{{\scriptsize $1$}}}
\put(133.00,17.00){\makebox(0,0)[cc]{{\scriptsize 63}}}
\put(120.00,13.00){\line(1,0){17.00}}
\put(141.00,13.00){\circle{8.00}}
\put(141.00,13.00){\makebox(0,0)[cc]{{$1$}}}
\put(141.00,3.00){\makebox(0,0)[cc]{$Sp_6(2)$}}
\put(76.00,13.00){\line(1,0){17.00}}
\end{picture}

\medskip

\end{itemize}

\item[{\rm (iii)}] $^4\Sg$ is the commuting graph of
$3$-transpositions in the Fischer group $Fi_{22}$ with the
suborbit diagram

\unitlength 0.85mm
\linethickness{0.4pt}
\begin{picture}(150.00,19.00)
\put(41.00,12.00){\circle{8.00}}
\put(62.00,12.00){\line(1,0){17.00}}
\put(98.00,12.00){\line(1,0){17.00}}
\put(41.00,12.00){\makebox(0,0)[cc]{{$1$}}}
\put(51.00,16.00){\makebox(0,0)[cc]{{\scriptsize $693$}}}
\put(76.00,16.00){\makebox(0,0)[cc]{{\scriptsize 1}}}
\put(104.00,16.00){\makebox(0,0)[cc]{{\scriptsize $512$}}}
\put(127.00,16.00){\makebox(0,0)[cc]{{\scriptsize $126$}}}
\put(88.50,12.00){\oval(19.00,8.00)[]}
\put(88.00,12.00){\makebox(0,0)[cc]{ $693$}}
\put(88.00,19.00){\makebox(0,0)[cc]{\scriptsize $180$}}
\put(141.00,12.00){\oval(18.00,8.00)[]}
\put(141.00,12.00){\makebox(0,0)[cc]{ $2816$}}
\put(141.00,19.00){\makebox(0,0)[cc]{\scriptsize $567$}}
\put(45.00,12.00){\line(1,0){17.00}}
\put(115.00,12.00){\line(1,0){17.00}}
\put(41.00,2.00){\makebox(0,0)[cc]{$2 \cdot U_6(2)$}}
\put(88.00,2.00){\makebox(0,0)[cc]{$2^{2+8}.U_4(2)$}}
\put(141.00,2.00){\makebox(0,0)[cc]{$U_4(3):2$}}
\end{picture}
\end{itemize}
\end{lem}

\pf By (\ref{21} (i)) $\Sg$ is locally $\Upsilon$. A connected graph
which is locally the complete multipartite graph $K_{t \times m}$ with
$t \ge 2$ parts of size $m \ge 2$ each, is isomorphic to $K_{(t+1)
\times m}$ (cf. Proposition 1.1.5 in
\cite{BCN89}) and hence (i) follows. For $t=2$ and $4$ let $\cP$ be the
polar space with diagram
$$\node_t\arc\node_t\darc\node_2$$
as in (\ref{20} (iii) (b), (c)) so that $\Upsilon$ is the collinearity
graph of $\cP$. Let $\cR$ be the geometry of rank 4, whose elements of
type 4 are the maximal complete subgraphs (of size $t^2+t+2$) in
$\Sg$; the elements of type 1 are the vertices, the
elements of type 2 are the edges and the elements of type 3 are the
complete subgraphs of size $t+2$ contained in more than one maximal
complete subgraph; the incidence relation is via inclusion. It is a standard fact
that $\cR$ is a $c$-extension of $\cP$ with the diagram
$$\node_1\stroke{\rm c}\node_t\arc\node_t\darc\node_2$$
If $t=2$ then by \cite{CP92} $\cR$ is a standard quotient of an
affinization of the polar space of $Sp_8(2)$. It is well known that
$\cP$ contains three classes of hyperplanes and we obtain the
possibilities (a), (b) and (c) in (ii) (notice that the antipodal
quotient of $^2\Sg^c$ is the complete graph which is not locally
$^2\Upsilon$). Finally the isomorphism type of $^4\Sg$ follows from
\cite{Pase94} and \cite{Pase95}.\qed

\medskip

Let us match the possibilities in (\ref{22}) with the examples we know.
Let $\Gm=\Gm(\cE)$ where $\cE$ is one of the seven examples of
$c.F_4(t)$-geometries in Section~\ref{s4}. If $\cE=\cE(3 \cdot
Fi_{22})$ then $\Gm(x,y) \cong $$^1\Sg$ and if $\cE=\cE(Fi_{22})$ then
$\Gm(x,y)$ consists of three connected components, each isomorphic to
$^1\Sg$. If $\cE=\cE(3 \cdot $$^2E_6(2))$ then $\Gm(x,y) \cong $$^2\Sg$$^a$
and if $\cE=\cE($$^2E_6(2))$ then $\Gm(x,y)$ consists of three connected
components, each isomorphic to $^2\Sg^a$. If $\cE=\cE(E_6(2))$ then
$\Gm(x,y) \cong ^2\Sg^b$ while if $\cE=\cE(2^{26}.F_4(2))$ then $\Gm(x,y)
\cong $$^2\Sg^c$. Finally if $\cE=\cE(B)$ then $\Gm(x,y) \cong
$$^4\Sg$.

\medskip

Let $\Sg$ be $^1\Sg$, $^2\Sg^a$, $^2\Sg^b$ or $^4\Sg$ as in (\ref{22})
(not $^2\Sg^c$) and $T$ be isomorphic to $S_3 \times S_3 \times S_3
\times S_3$, $\Omega_8^-(2):2$, $\Omega_8^+(2):2$ or $Fi_{22}$,
respectively, so that $\Sg$ is isomorphic to the commuting graph of a conjugacy class of
$3$-transpositions in $T$. This means that there is a bijection
$$j: \Sg \to T,$$
such that $o(j(u)j(v))=2$ if $v \in \Sg(u)$ and
$o(j(u)j(v))=3$ if $v \in \Sg_2(u)$. Notice that by
(\ref{21} (ii)) $u$ and $v$ (considered as involutions in
$F$) satisfy $o(uv)=2$ if $v \in \Sg(u)$ and $o(uv)=3$ if $v \in
\Sg_2(u)$. For $t \in T$ and $x \in \Sg$ the action $x^t=j^{-1}(j(x)^t)$
turns $T$ into an automorphism group of $\Sg$.

\begin{lem} \label{23}
In the above terms (i.e. with $\Sg \not \cong $$^2\Sg^c$)
\begin{itemize}
\item[{\rm (i)}] whenever $u,v \in
\Sg$ and $w=uvu$, we have $j(w)=j(u)j(v)j(u)$;
\item[{\rm (ii)}] if $u \in \Sg$, then $u$ (considered as an
automorphism of $\Dl$) stabilizes $\Sg$ as a whole.
\end{itemize}
\end{lem}

\pf Let us proceed with (i) (which will immediately imply (ii)).
If $u$ and $v$ are adjacent then $w=v$,
$j(w)=j(v)$ and the assertion is obvious. If $v \in
\Sg_2(u)$ we define $x$ so that $j(x)=j(u)j(v)j(u)$ (this is possible
since $j(\Sg)$ is a conjugacy class in $T$) and all we have to
show is that $w:=uvu$ coincides with $x$. Let $\Sg(u,v)$ be the
$\mu$-graph in $\Sg$ (the subgraph induced by $\Sg_1(u) \cap \Sg_1(v)$).
Since $j(u)$ fixes $\Sg_1(u)$ elementwise and maps $v$ onto $x$, we
observe that $\Sg(u,x)=\Sg(u,v)$. Furthermore,
$j(\Sg(u,v))=j(\Sg) \cap T(u) \cap T(v)$. In
view of the suborbit diagrams in (\ref{22}) this shows, particularly,
that $\Sg(u,v)$ contains a pair of non-adjacent vertices, say $p$ and
$q$. Then $\{u,v,x\} \subset \Sg_1(p) \cap \Sg_1(q)$ and $q \in \Sg_2(p)$.
On the one hand $\Sg(p,q)$ is isomorphic to the commuting graph of
involutions in $j(\Sg) \cap (T(p) \cap T(q))$ and on the other hand
$\Sg(p,q)$ is a subgraph in the commuting graph $\Dl(p,q)$ of
involutions in $F(p) \cap F(q)$. Thus the restriction $\sigma$ of $j^{-1}$ to
$\Sg(p,q)$ induces an isomorphic embedding of the commuting graph
$\Pi$ of involutions in $j(\Sg(p,q))$ into the commuting graph
$\Upsilon$ (as in (\ref{20})) of involutions in $\Dl(p,q)$ and
it remains to show that the embedding $\sigma$ is {\it $S_3$-consistent}
in the sense that it maps the set
$\{u,v,x\}$ of involutions in an $S_3$-subgroup in
$K:=T(p) \cap T(q)$ onto the set of involutions in an $S_3$-subgroup
in $H=F(p) \cap F(q)$. If $t=1$ then $|\Pi|=|\Upsilon|=9$,
$\sigma$ is an isomorphism and a triple of vertices in
$\Pi$ or in $\Upsilon$
generate an $S_3$-subgroup if and only if these three vertices are
pairwise non-adjacent, the result is trivial in this case.

In the case $t=2$ the graph $\Pi$ is the commuting graph of
$3$-transpositions in $K \cong \Omega_6^+(2):2 \cong S_8$ with the suborbit
diagram

\medskip

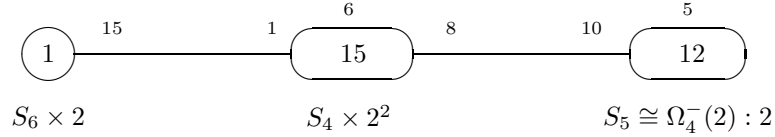
\unitlength 0.85mm
\linethickness{0.4pt}
\begin{picture}(148.00,19.00)
\put(39.00,12.00){\circle{8.00}}
\put(60.00,12.00){\line(1,0){17.00}}
\put(96.00,12.00){\line(1,0){17.00}}
\put(39.00,12.00){\makebox(0,0)[cc]{{$1$}}}
\put(49.00,16.00){\makebox(0,0)[cc]{{\scriptsize $15$}}}
\put(74.00,16.00){\makebox(0,0)[cc]{{\scriptsize 1}}}
\put(102.00,16.00){\makebox(0,0)[cc]{{\scriptsize $8$}}}
\put(124.00,16.00){\makebox(0,0)[cc]{{\scriptsize $10$}}}
\put(86.50,12.00){\oval(19.00,8.00)[]}
\put(86.00,12.00){\makebox(0,0)[cc]{ $15$}}
\put(86.00,19.00){\makebox(0,0)[cc]{\scriptsize $6$}}
\put(139.00,12.00){\oval(18.00,8.00)[]}
\put(139.00,12.00){\makebox(0,0)[cc]{ $12$}}
\put(43.00,12.00){\line(1,0){17.00}}
\put(113.00,12.00){\line(1,0){17.00}}
\put(39.00,2.00){\makebox(0,0)[cc]{$S_6 \times 2$}}
\put(86.00,2.00){\makebox(0,0)[cc]{$S_4 \times 2^2$}}
\put(139.00,2.00){\makebox(0,0)[cc]{$S_5\cong \Omega_4^-(2):2$}}
\put(139.00,19.00){\makebox(0,0)[cc]{\scriptsize $5$}}
\end{picture}

\medskip

\noindent
or the commuting graph of $3$-transpositions in $K \cong
\Omega_6^-(2):2 \cong U_4(2):2$ with the suborbit diagram

\medskip

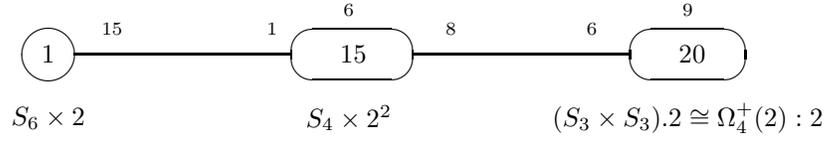
\unitlength 0.85mm
\linethickness{0.4pt}
\begin{picture}(148.00,19.00)
\put(39.00,12.00){\circle{8.00}}
\put(60.00,12.00){\line(1,0){17.00}}
\put(96.00,12.00){\line(1,0){17.00}}
\put(39.00,12.00){\makebox(0,0)[cc]{{$1$}}}
\put(49.00,16.00){\makebox(0,0)[cc]{{\scriptsize $15$}}}
\put(74.00,16.00){\makebox(0,0)[cc]{{\scriptsize 1}}}
\put(102.00,16.00){\makebox(0,0)[cc]{{\scriptsize $8$}}}
\put(124.00,16.00){\makebox(0,0)[cc]{{\scriptsize $6$}}}
\put(86.50,12.00){\oval(19.00,8.00)[]}
\put(86.00,12.00){\makebox(0,0)[cc]{ $15$}}
\put(86.00,19.00){\makebox(0,0)[cc]{\scriptsize $6$}}
\put(139.00,12.00){\oval(18.00,8.00)[]}
\put(139.00,12.00){\makebox(0,0)[cc]{ $20$}}
\put(43.00,12.00){\line(1,0){17.00}}
\put(113.00,12.00){\line(1,0){17.00}}
\put(39.00,2.00){\makebox(0,0)[cc]{$S_6 \times 2$}}
\put(86.00,2.00){\makebox(0,0)[cc]{$S_4 \times 2^2$}}
\put(139.00,2.00){\makebox(0,0)[cc]{$(S_3 \times S_3).2\cong \Omega^+_4(2):2$}}
\put(139.00,19.00){\makebox(0,0)[cc]{\scriptsize $9$}}
\end{picture}

\medskip

By (\ref{20} (iii) (b)) $\Upsilon$ is the
collinearity graph of the polar space $\cP=\cP(Sp_6(2))$. The lines of this
space can be identified with a class of triangles in $\Upsilon$ so that
every edge is contained in a unique line and a vertex outside a line is
adjacent to 1 or 3 vertices on the line. It is easy to see that for
every triangle in $\Pi$ there is a vertex adjacent to 2 vertices in the
triangle. Hence $\sigma(\Pi)$ does not contain lines of $\cP$
and since the valency of $\Pi$ (which is 15) equals to the number of
lines containing a given point, every line intersects
$\sigma(\Pi)$ in 2 or 0 vertices. Hence $\sigma(\Pi)$ is the complement
to a hyperplane in $\cP$. All hyperplanes in $\cP$ are known, in
particular every hyperplane with complement of size 28, respectively 36
has stabilizer in $H \cong Sp_6(2)$ isomorphic to $\Omega_6^+(2):2$ or
$\Omega_6^-(2):2$, which shows that $\sigma$ is $S_3$-consistent.

\medskip

In the case $t=4$ $\Pi$ is the commuting graph of $3$-transpositions
in $K \cong U_4(3):2$ with the suborbit diagram

\medskip

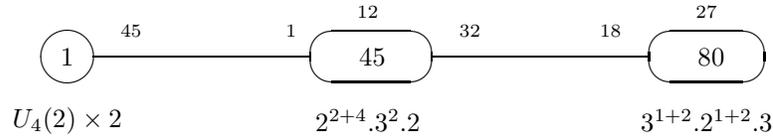
\unitlength 0.85mm
\linethickness{0.4pt}
\begin{picture}(148.00,19.00)
\put(39.00,12.00){\circle{8.00}}
\put(60.00,12.00){\line(1,0){17.00}}
\put(96.00,12.00){\line(1,0){17.00}}
\put(39.00,12.00){\makebox(0,0)[cc]{{$1$}}}
\put(49.00,16.00){\makebox(0,0)[cc]{{\scriptsize $45$}}}
\put(74.00,16.00){\makebox(0,0)[cc]{{\scriptsize 1}}}
\put(102.00,16.00){\makebox(0,0)[cc]{{\scriptsize $32$}}}
\put(124.00,16.00){\makebox(0,0)[cc]{{\scriptsize $18$}}}
\put(86.50,12.00){\oval(19.00,8.00)[]}
\put(86.00,12.00){\makebox(0,0)[cc]{ $45$}}
\put(86.00,19.00){\makebox(0,0)[cc]{\scriptsize $12$}}
\put(139.00,12.00){\oval(18.00,8.00)[]}
\put(139.00,12.00){\makebox(0,0)[cc]{ $80$}}
\put(43.00,12.00){\line(1,0){17.00}}
\put(113.00,12.00){\line(1,0){17.00}}
\put(39.00,2.00){\makebox(0,0)[cc]{$U_4(2) \times 2$}}
\put(86.00,2.00){\makebox(0,0)[cc]{$2^{2+4}.3^2.2$}}
\put(139.00,2.00){\makebox(0,0)[cc]{$3^{1+2}.2^{1+2}.3$}}
\put(139.00,19.00){\makebox(0,0)[cc]{\scriptsize $27$}}
\end{picture}

\medskip

\noindent
and locally it is the point graph of the generalized quadrangle of
order $(4,2)$ associated with $U_4(2)$. This shows that the vertices,
edges and maximal cliques (of size 6) with respect to inclusion form a
geometry $\cG$ with the diagram $\node_1\stroke{\rm
c}\node_4\darc\node_2$. Recall that in the considered situation
$\Upsilon$ is the collinearity graph of the polar space $\cP=\cP(U_6(2))$
so that the lines and planes in $\cP$ can be identified with certain 
complete subgraphs of size 5 and the maximal complete subgraphs (of
size 21) in $\Upsilon$. We extend $\sigma$ to a morphism
$\delta$ of $\cG$ into $\cP$ as
follows. For $t \in \Pi$ put $\delta(t)=\sigma(t)$; for an edge $e$ of $\Pi$ define $\delta(e)$ to be the unique
line in $\cP$ containing $\sigma(e)$ and for a 6-clique $X$ in $\Pi$
define $\delta(X)$ to be the unique plane in $\cP$ containing
$\sigma(X)$. In the polar space $\cP$ a point outside a line
is adjacent to 1 or 5
points of the line, while in $\Pi$ for every triangle there is a vertex
adjacent to exactly 2 vertices in the triangle. It is an easy
combinatorial exercise to check, using this observation, that
$\delta$ is an isomorphic embedding. The distribution
diagram of the graph on the elements of type 3 in $\cG$ in which two
such elements are adjacent if they are incident to a common element of
type 2 (cf. p. 202 in \cite{I99}) shows, that every edge in this graph is
contained in a unique triangle. Hence the image under $\delta$ of the
set of elements of type 3 from $\cG$ is a hyperplane in the dual polar
space associated with $\cP$. Now by Proposition 2.1
in \cite{C94} the stabilizer of $\delta(\cG)$ in $H \cong U_6(2):S_3$
is isomorphic to $K \cong U_4(3):2$ and hence $\sigma$ is
$S_3$-consistent.\qed

\medskip

Let us try to identify $\Sg$ as a subgraph in $\Dl$. We proceed to
this via identifying the stabilizer $S$ of $\Sg$ in $F$. By the
definition, $S$ is the normalizer in $F$ of $\Sg$ where the latter
is considered as a set of central involutions in $F$. The following result
in a consequence of (\ref{23} (ii)).

\begin{cor} \label{23.5}
Let $U$ be the subgroup in $F$ generated by the involutions from $\Sg$.
If $\Sg \not \cong $$^2\Sg^c$ then $S=N_F(U)$ and $U$ is isomorphic to the
group $T$ as in the proof of $(\ref{23})$.\qed
\end{cor}

In order to identify $S$ precisely we treat different values of $t$
separately.

\medskip

In the case $t=1$ let $p \in \Dl$, $q \in \Psi_3(p)$ and $r$ be the
image of $q$ under conjugation by $p$. Then by the proof of
(\ref{23}), up to conjugation in $F$ we have
$^1\Sg=\{p,q,r\} \cup (\Dl(p) \cap \Dl(q))$. This specifies $^1\Sg$ and
in view of the list of maximal subgroups in $^1F$ gives the following.

\begin{lem} \label{24}
The stabilizer $^1S$ of $^1\Sg$ in $^1F \cong \Omega_8^+(2).S_3$ is
isomorphic to $S_3 \wr S_4$. Furthermore $^1S$ is a subgroup of index
$3$ in a maximal subgroup of $^1F$.
\end{lem}

Let us turn to the case $t=2$. Recall that $^2F \cong F_4(2)$ contains
two conjugacy classes of subgroups isomorphic to $Sp_8(2)$. The
representatives $P_1$ and $P_2$ of these classes can be chosen in such way
that the centralizer of a central involution in $P_1$ (isomorphic to
$2^7.Sp_6(2)$) stabilizes a point in $\cF$ while the similar
centralizer in $P_2$ stabilizes a symplecton. Also $^2F$ has two
irreducible 26-dimensional $GF(2)$-modules $V_1$ and $V_2$ such that
$P_i$ fixes a vector in $V_i$ and acts irreducibly on $V_{3-i}$. The
classes containing $P_1$ and $P_2$ as well as the modules $V_1$ and
$V_2$ are permuted by the outer (diagram) automorphism of $^2F$. Under
the choice of the classes, the class $M$ of 255 central
involutions in $P_1$ (which is a class of $3$-transpositions) 
is a subset of $\Dl$ and the subgraph induced by $M$ is clearly the
commuting graph. For $\varepsilon=+$ or $-$ let $Q^\varepsilon$ be
(the unique up to conjugation) subgroup in $P_1$ isomorphic to
$\Omega^{\varepsilon}_8(2):2$. Then the involutions from $M$
contained in $Q^\varepsilon$ induce a subgraph isomorphic to $^2\Sg^a$
or $^2\Sg^b$ if $\varepsilon=+$ or $-$, respectively.

\begin{lem} \label{25}
For $\alpha=a$ or $b$ the stabilizer $^2S^\alpha$ of $^2\Sg^\alpha$ in
$^2F$ is a conjugate of $Q^\varepsilon$ for $\varepsilon=+$ or $-$,
respectively.
\end{lem}

\pf By Theorem 4 in \cite{LS98} every subgroup in $F_4(2)$ isomorphic to
$\Omega_8^{\pm}(2):2$ stabilizes a vector either in $V_1$ or in $V_2$
and by  Table 2 in \cite{CC88} such a
subgroup is contained either in a conjugate of $P_1$ or in a conjugate
of $P_2$.\qed

\medskip 

Notice that $N_{^2F}(O^2(Q^+)) \cong \Omega_8^+(2):S_3$ and
that for $u \in M$ the subgraph in $\Dl$ induced by the involutions
from $M$ which do not commute with $u$ is isomorphic to $^2\Sg^c$.

\medskip

Finally consider the case $t=4$. It has been shown by B.~Fischer (cf.
\cite{Coo83}) that $^2E_6(2)$ contains at least three conjugacy classes
of subgroups isomorphic to $Fi_{22}$ and that these three classes are
fused in $^4F \cong $$^2E_6(2).S_3$. By \cite{N92} the subgroups in $^4F$
isomorphic to $Fi_{22}$ are all conjugate.

\begin{lem} \label{26}
All subgroups in $^4F$ isomorphic to $F_{22}$ are conjugate and the
stabilizer $^4S$ of $^4\Sg$ is the normalizer of such a subgroup,
isomorphic to $Fi_{22}:2$.\qed
\end{lem}

Notice that $^4S$ is contained in a subgroup of index 3 in $^4F$.

\medskip

We have the following result (cf. \cite{IS96} and Lemma 5.10.8 in
\cite{I99}) 

\begin{lem} \label{A} 
Let $S=$$^4S\cong Fi_{22}:2$ be the stabilizer of $\Sg=$$^4\Sigma$ in $F=$$^4F \cong
$$^2E_6(2):S_3$. Then 
\begin{itemize}
\item[{\rm (i)}] $S$, acting on $\Dl$ has four orbits:
$\Omega_1=\Sigma$, $\Omega_2$, $\Omega_3$ and $\Omega_4$; 
\item[{\rm (ii)}]
$|\Omega_i|=3\,510,~~142\,155,~~694\,980~~and~~3\,127\,410,$ for
$i=1,2,3$ and $4$, respectively; 
\item[{\rm (iii)}] if $w_i \in \Omega_i$, then $S(w_i)$ is isomorphic, 
to $2 \cdot
U_6(2).2,~~2^{10}.M_{22}.2,~~2^7:Sp_6(2)~~and~~2.2^9:P\Sigma L_3(4)$
for $i=1,2,3$ and $4$, respectively; 
\item[{\rm (iv)}] for $i=2,3$ and $4$ the vertex $w_i$ is
adjacent in $\Delta$, respectively, to $22$, $126$ and $22$ vertices
from $\Sigma$;  
\item[{\rm (v)}] $S$ acts transitively on the set of 
maximal cliques $($of size $22)$ in $\Sigma$; $\Delta(w_2) \cap \Sigma$
and $\Delta(w_4) \cap \Sigma$ are such cliques and $S(w_2) \cong 2^{10}.M_{22}.2$ 
is the stabilizer of $\Delta(w_2) \cap \Sg$ in $S$; 
\item[{\rm (vi)}] $\Dl(w_2) \cap \Sg \subset O_2(F(w_2)) \setminus \{w_2\}$ 
while
$|\Dl(w_4) \cap \Sg \cap O_2(F(w_4))|=1$.\qed 
\end{itemize}
\end{lem}

Notice that $O_2(S(w_2))$ is the 10-dimensional Todd module for
$S(w_2)/O_2(S(w_2))$.

\begin{lem} \label{AB}
In terms of $(\ref{A})$ the subgroup $L=\la \Dl(w_2) \cap \Sg \ra=O_2(S(w_2))$ 
is contained in 
exactly two conjugates of $S$ in $F$, namely, in $S$ and in $w_2Sw_2$,
moreover, $S \cap w_2Sw_2 =S(w_2) \cong 2^{10}.M_{22}.2$.
\end{lem} 

\pf By (\ref{A} (v)) $N_F(L)$ acts transitively on the set of conjugates of
$S$ in $F$, containing $L$. By (\ref{A} (vi)) $w_2$ is the unique vertex of
$\Dl$ such that $L \le O_2(F(w_2))$ and hence $N_F(L) \le F(w_2)$. Since
$O_2(F(w_2))$ is extraspecial of order $2^{21}$,
$N_{O_2(F(w_2))}(L)=L\la w_2\ra$. We have
$$N_S(L)O_2(F(w_2))/O_2(F(w_2)) \cong M_{22}.2$$
and by the list of
maximal subgroups in $U_6(2)$ \cite{ATLAS} the inclusion $M_{22}.2 < N
\le U_6(2).S_3$ implies $N \ge U_6(2).2$. Hence 
$N_{F}(L)=N_{S}(L)\la w_2 \ra$ and the result follows.\qed 

\medskip

The next proposition is an immediate consequence of (\ref{A}),
(\ref{13.5}) and (\ref{22} (iii)). 

\begin{prop} \label{B}
Suppose that $t=4$ and let $y \in \Gm_2^3(x)$. Then
\begin{itemize}
\item[{\rm (i)}] the subgraph $\Gm(x,y)$ is connected with $3\,510$ vertices;
\item[{\rm (ii)}] if $u \in \Gm(y) \setminus \Gm(x,y)$, then $u \in
\Gm_2(x)$ and $u$ is adjacent in $\Gm$ to a vertex in $\Gm(x,y)$;
\item[{\rm (iii)}] $\Gm_2^3(x) \cap \Gm_2^4(x)=\emptyset$;
\item[{\rm (iv)}] $|\Gm_2^3(x)|=|\Dl|\cdot |\Psi_3(x)|/3\,510$.\qed
\end{itemize}
\end{prop}

The following lemma provides some further details on (\ref{19} (ii)). 

\begin{lem} \label{C}
Let $e$, $f$ be distinct elements of type $2$ in $\cE$, such that
$\Pi:=\Xi^e \cap \Xi^f \ne \emptyset$. Let $\widetilde \Pi$ be the
image of $\Pi$ under the natural mapping of $\Xi^e$ onto $\widetilde
\Xi^e$. Assume without loss that $e=\{x,y\}$, $f=\{x,z\}$. Then the
following assertions hold:

\begin{itemize}
\item[{\rm (i)}] if $i_x(z) \in \Psi_1(i_x(y))$ then $|\Pi|=2 \cdot
|\widetilde \Pi|=4(t^2+t+2)$ and $z \in \Xi^e$;
\item[{\rm (ii)}] if $i_x(z) \in \Psi_2^2(i_x(y))$ then $|\Pi|=|\widetilde
\Pi|=3(2t+1)+1$ and $|\Gm(z) \cap \Xi^e|=6(2t+1)+2$;
\item[{\rm (iii)}] if $i_x(z) \in \Psi_2^4(i_x(y))$ then $|\Pi|=|\widetilde
\Pi|=2$ and $|\Gm(z) \cap \Xi^e|=2(t^2+t+1)+1$;
\item[{\rm (iv)}] if $i_x(z) \in \Psi_3(i_x(y))$ then $|\Pi|=|\widetilde
\Pi|=1$ and $\Gm(z) \cap \Gm(x) \cap \Xi^e=\emptyset$; if in addition
$t=4$, then $\Gm(z) \cap \Xi^e=\{x\}$.
\end{itemize}

\noindent
In particular if $t=4$ then $\Xi^e$ is geodetically closed.

\end{lem}

\pf By (\ref{17}) and (\ref{B}) if $z$ is adjacent to a vertex from
$\Xi^e \setminus (\{x\} \cup \Gm(x))$, then either $z \in \Xi^e$, or 
$t \ne 4$ and $i_x(z) \in \Psi_3(i_x(y))$. Now the result can be easily 
deduced from Fig. 2, Fig. 3, the fact that 
$$\Gm(x) \cap \Xi^e=\{y\} \cup \{u \mid i_x(u) \in \Psi_1(i_x(y))\}$$
and the definition of $\widetilde \Xi^e$.\qed 

\medskip

Notice that if $t=4$, $^4G \cong Co_2$ is the automorphism group of
$\widetilde \Xi^e$, $\widetilde \Pi$ is as in (\ref{C} (i)) and
if $M$ is the stabilizer of $\widetilde \Pi$ in $^4G$, then $M \cong
2^{10}.M_{22}.2$, where $O_2(M)$ is the Golay code module for
$M/O_2(M)$.  

\section{A characterization of the Baby Monster graph} \label{s8}
In this section (which consists of a few subsections) we prove
Theorems~\ref{aaa} and \ref{bbb} by showing that $\Gm($$^4\cE) \cong
\Th$ and that $^4\cE \cong \cE(B)$.
Thoughout the section we assume
that $t=4$. 

\subsection{$x$-equivalence on $\Gm_2(x)$}
Let $x \in \Gm$, $y \in \Gm^\alpha_2(x)$ for $\alpha=3$ or $4$, $\Gm(x,y)$ be the corresponding
$\mu$-graph, $\Sg^y=\{i_x(u) \mid u \in \Gm(x,y)\}$ and $S^y$ be the
stabilizer of $\Sg^y$ in $F \cong $$^2E_6(2).S_3$. Then by (\ref{14}),
(\ref{16}), (\ref{26}), (\ref{A}), $S^y$ is specified by $\alpha$ up to
conjugation in $F$ and $S^y$ determines $\Sg^y$ uniquely. We summarize
this in the following proposition.

\begin{prop} \label{D}
The following assertions hold:
\begin{itemize}
\item[{\rm (i)}] if $y \in \Gm^4_2(x)$, then $S^y \cong
2^{1+20}_+:U_4(3).2^2$ is contained in $F(p) \cong
2^{1+20}_+:U_6(2).S_3$ for a unique point $p$ of $\cF$ and $\Sg^y$ is
the unique orbit of length $648$ of $S^y$ on $\Dl$ (contained in
$\Psi_1(p)$); 
\item[{\rm (ii)}] if $y \in \Gm_2^3(x)$, then $S^y \cong Fi_{22}:2$ and
$\Sg^y$ is the unique orbit of length $3\,510$ of $S^y$ on $\Dl$; 
\item[{\rm (iii)}] in either case if $H^y=S^yF^\infty$, then 
$H^y \cong $$^2E_6(2):2$ is a subgroup of index $3$ in $F$.\qed
\end{itemize}
\end{prop}

In the next subsection we show that $H^y$ is independent on the
choice of $y \in \Gm_2(x)$. We need another preliminary result.

\begin{lem} \label{DE}
Let $F(p) \cong 2^{1+20}_+:U_6(2).S_3$ be the stabilizer in $F$ of a
point $p$ of $\cF$, $Z=\la p\ra$, $Q=O_2(F(p))$, $I=F(p)^{\infty}$.
Then $\bar I:=I/Z \cong 2^{20}:U_6(2)$ has four classes $\cU_i$, $0 \le
i \le 4$, of complements to
$Q/Z$, such that $F(p)/I \cong S_3$ normalizes $\cU_0$ and permutes transitively
the remaining classes; the preimage in $I$ of a subgroup from $\cU_i$
splits over $Z$ if and only if $i=0$.\qed
\end{lem}

Let  $p \in F$, $M_1$, $M_2$, $M_3$ be representatives of the
conjugacy classes of $Fi_{22}$-subgroups in $F^\infty \cong $$^2E_6(2)$
(compare the paragraph before (\ref{26})) such that $p \in M_i$ for $1
\le i \le 3$ and $q \in \Psi_3(p)$. Then (within a suitable ordering)
we have $(F(p) \cap F(q))^\infty \in \cU_0$, $M_i(p)/\la p \ra \in \cU_i$ for
$1 \le i \le 3$. 

\begin{lem} \label{E}
For every $y \in \Gm_2(x)$ there is exactly one vertex $y^\prime \in \Gm_2(x)
\setminus \{y\}$ such that $\Sg^y=\Sg^{y^\prime}$. Furthermore, 
$i_z(y')=i_z(x)i_z(y)i_z(x)$ for every $z \in \Gm(x,y)$. 
\end{lem}

\pf Suppose first that $y \in \Gm_2^4(x)$. Let $e=\{x,z\}$ be the
element of type 2 in $\cE$ such that $y \in \Xi^e$ and $\chi:\Xi^e \to
\widetilde \Xi^e$ be the natural mapping. Since $\Xi^e$ is the 2-clique
extension of $\widetilde \Xi^e$, it is clear that $\Sg^y=\Sg^{y^\prime}$
whenever $\chi(y)=\chi(y^\prime)$. On the other hand, it is easy to see from Fig.
2 and Fig. 3 that $z$ is the only vertex in $\Gm(x)$ adjacent to every
vertex in $\Gm(x,y)$. Hence $\Sg^{y^\prime}=\Sg^y$ implies $y^\prime \in \Xi^e$.
Finally, different pairs of vertices at distance 2 in $\widetilde
\Xi^e$ have different $\mu$-graphs, which gives $\chi(y^\prime)=\chi(y)$. 

Now suppose that $y,y^\prime \in \Gm_2^3(x)$, $y \ne y^\prime$, $\Sg^y=\Sg^{y^\prime}$
and $z \in \Gm(x,y)$. Put $p=i_z(x)$, $q=i_z(y)$, $q^\prime=i_z(y^\prime)$. We claim
that $q^\prime=pqp$. The subgroup $U:=O^{\infty}(F(p) \cap F(q)) \cong
U_6(2)$ is generated by the involutions in $\Dl(p) \cap \Dl(q)=\Dl(p)
\cap \Dl(q^\prime)$. Hence $U$ commutes with the unique element $r \in
O_2(F(p))$ which maps $q$ onto $q^\prime$. Since $U$ acts irreducibly
on $O_2(F(p))/\la p \ra$ and $q\ne q^\prime$, we have $r=p$, hence 
the claim follows and proves the uniqueness statement in the lemma. Now let $y$ and $z$ be as
above, and let $y^\prime$ be such that 
$i_z(y^\prime)=i_z(x)i_z(y)i_z(x)$.
We claim that $S^y=S^{y^\prime}$ (then (\ref{D} (ii)) will 
imply the equality $\Sg^y=\Sg^{y^\prime}$). The subgroup $K$
of $F$ generated by 
$$\{i_x(u) \mid u \in \{z\} \cup (\Gm(x) \cap \Gm(y)
\cap \Gm(z)) \}$$
is isomorphic to $2 \cdot U_6(2)$ (a non-split
extension) and $K=O^2(S^y(r))=O^2(S^{y^\prime}(r))$, where $r=i_x(z)$.
By (7.8) there is $f \in F$ which conjugates $S^y$ onto $S^{y^\prime}$.
Since $S^y$ acts transitively on the set of vertices of $\Sg^y$ and $r$
is one of these vertices, we can choose $f$ to normalize $K$. By (\ref{A})
$r$ is the unique point in $\cF$, stabilized by $K$. This implies that
$N_{F}(K) \le F(r)$. Since $K$ does not split over its centre, by
(\ref{DE}) $N_{F(r)}(K) \cong 2 \cdot U_6(2).2 \cong N_{S^y}(K)$. Hence
$S^{y^\prime}=S^y$.\qed 
 
\medskip

As a consequence of the proof of (\ref{E}) and in view of the fact
that different pairs of vertices at distance 2 in 
$\widetilde \Xi^e$ have different $\mu$-graphs, we obtain the following. 

\begin{lem} \label{EF}
Suppose that $y_1,y_2 \in \Gm_2(x)$ and for some $z \in \Gm(x,y_1) \cap
\Gm(x,y_2)$ we have $\Gm(x) \cap \Gm(y_1) \cap \Gm(z)=\Gm(x) \cap
\Gm(y_2) \cap \Gm(z)$. Then $\Gm(x,y_1)=\Gm(x,y_2)$.\qed
\end{lem}

As a direct consequence of (\ref{EF}) we obtain the following result
(which can also be checked in $\Dl$ directly).

\begin{lem} \label{8.5} 
Let $p \in \Dl$, $q,q^\prime \in \Dl_2(p)$ and suppose that $\Dl(p)
\cap \Dl(q)=\Dl(p) \cap \Dl(q^\prime)$. Then either $q^\prime=q$ or
$q^\prime=pqp$.\qed
\end{lem}

To the end of this subsection we introduce some important notions.

\begin{defi} \label{defi}
Let $x \in \Gm$. Two vertices 
$y,y^\prime \in \Gm_2(x)$ are said to be $x$-equivalent if
$\Gm(x,y)=\Gm(x,y^\prime)$. By $\pi_x$ we denote the permutation of
the vertex-set of $\Gm$ which fixes every vertex in $\Gm \setminus
\Gm_2(x)$ and swaps the pairs of $x$-equivalent vertices in $\Gm_2(x)$.
\end{defi}

By (\ref{E}) each class of $x$-equivalent vertices is of size 2, so
that $\pi_x$ is well defined. In the Baby Monster
graph $\Theta$ two vertices $b,b^\prime \in \Th_2(a)$ are $a$-equivalent 
if $\la a,b\ra=\la a,b^\prime \ra$ (which happens exactly when $b^\prime=aba$.) 
Furthermore, $a,b \in
\Theta$ do not commute if and only if $b \in \Theta_2(a)$, so that
$\pi_a$ is nothing, but the conjugation by $a$, and clearly it is an
automorphism of $\Theta$. 

\medskip

Our ultimate goal is to prove the following.

\begin{prop} \label{ult}
For every $x \in \Gm$ the permutation $\pi_x$ is an automorphism of $\Gm$.
\end{prop}

The proof of (\ref{ult}) will be achieved by showing that
restrictions of $\pi_x$ to various subsets of vertices of $\Gm$ are
automorphisms of the subgraphs induced by that subsets. Our first
result of this type follows directly from (\ref{E}).

\begin{lem} \label{pi0}
If $z \in \Gm(x)$ and $y \in \Gm(z)$, then
$i_z(y^{\pi_x})=i_z(x)i_z(y)i_z(x)$, in particular the restriction 
$\pi_x^z$ of $\pi_x$ to $\Gm(z)$ is an automorphism of the subgraph in
$\Gm$ induced by $\Gm(z)$.\qed
\end{lem}

The permutation $\pi_x$ fixes every vertex in $\{x\} \cup \Gm(x)$ and
permutes the pairs of $x$-equivalent vertices in $\Gm_2(x)$. 
Hence by (\ref{pi0}) we also have the following. 

\begin{lem} \label{pi1}
For every $z \in \Gm(x)$ the restriction of $\pi_x$ to the set $\Gm(x) \cup
\Gm(z)$ is an automorphism of the subgraph induced by this set.\qed 
\end{lem}

\subsection{The class of $\mu$-graphs}
In this subsection we specify the family 
$$\cS=\{\Sg^u \mid u \in \Gm_2(x)\}$$
of subgraphs in $\Dl$ up to simultaneous conjugation 
by automorphisms of $\Dl$. 

\medskip

For $y \in \Gm^\alpha_2(x)$ let us adopt the notation introduced in
(\ref{D}) and in the paragraph before that lemma.  By (\ref{E}) there
are exactly $|\Gm_2^\alpha(x)|/2$ distinct $\mu$-graphs $\Gm(x,u)$ for
$u \in \Gm_2^\alpha(x)$ and by (\ref{19}),
(\ref{B}), (\ref{D}) and (\ref{E}) we have 

$$|\Gm_2^\alpha(x)|/2=[F:S^y]/3=[H^y:S^y].$$
Thus exactly one third of the images of $\Sg^y$ under the elements $f
\in F$ is contained in $\cS$ and we are going to show that such an 
image is contained in $\cS$ exactly when $f \in H^y$. First
of all by the above equality we obtain the following.

\begin{lem} \label{F}
For $\alpha=3$ or $4$ let $y \in\Gm_2^\alpha(x)$. Then the following
two statements are equivalent 

\begin{itemize}
\item[{\rm (i)}] $H^u=H^y$ for every $u \in \Gm_2^\alpha(x)$; 
\item[{\rm (ii)}] for every $h \in H^y$ the image of $\Sg^y$ under $h$ 
coincides with $\Sg^u$ for some $u \in \Gm_2^\alpha(x)$.\qed
\end{itemize}
\end{lem}

\begin{lem} \label{G}
The equivalent statements $(i)$ and $(ii)$ in $(\ref{F})$ hold.
\end{lem}

\pf Suppose first that $\alpha=3$. For a vertex $z \in \Gm(x,y)$ let
$p=i_z(x)$, $q=i_z(y)$ (so that $q \in \Psi_3(p)$), $r \in \Psi_3(p)
\cap \Psi_1(q)$ and $u=i^{-1}_z(r)$. Then $u \in \Gm_2^3(x) \cap
\Gm(y)$. By (2.1) $F(p) \cap F(q) \cap F(r) \cong 2^9.L_3(4).S_3$ is a
maximal parabolic in $F(p) \cap F(q) \cong U_6(2).S_3$. Hence $\Dl(p)
\cap \Dl(q) \cap \Dl(r)$ consists of 21 3-transpositions from the set
$\Dl(p) \cap \Dl(q)$ of 693 3-transpositions of $(F(p) \cap
F(q))^\infty \cong U_6(2)$. This means that $\Gm(x)
\cap \Gm(y) \cap \Gm(u)$ contains a set of 22 pairwise adjacent vertices, namely
$$\{z\} \cup i_z^{-1}(\Dl(p) \cap \Dl(q) \cap \Dl(r)).$$
On the one hand by (\ref{E}) $\Gm(x,y) \ne \Gm(x,u)$ (equivalently $S^y \ne
S^u$) and on the other hand by (\ref{A} (v)) and (\ref{AB}) $S^y \cap S^u \cong
2^{10}.M_{22}.2$ and $S^y$ is conjugate to $S^u$ in $F^\infty$, which
implies $H^u=H^y$, since $H^y=F^\infty S^y$. 
By the construction for every $22$-vertex complete
subgraph $\Upsilon$ in $\Gm(x,y)$ there is (a unique) $u \in \Gm_2^3(x)
\cap \Gm(y)$ such that $\Gm(x,y) \cap \Gm(u)=\Upsilon$. Since a
$Fi_{22}$-subgroup is maximal in $F^\infty \cong $$^2E_6(2)$, the graph
on the set of $F^\infty$-conjugates of $S^y$ in which two such
subgroups are adjacent if they share a set of 22 pairwise commuting
3-transpositions, is connected. This implies (\ref{F} (ii)) for $\alpha=3$.

\medskip

Now let $y \in \Gm_2^4(x)$, so that $y \in \Xi^e$ for some $e=\{x,z\}$
and $S^y \le F(p)$ where $p=i_x(z)$. The stabilizer of $e$ in the
automorphism group $Co_2$ of $\widetilde \Xi^e$ is isomorphic to $U_6(2):2$
(cf. (\ref{11} (iii)) and hence there is a unique subgroup $D^z$ of index $3$
in $F$ such that for every $d \in D^z(p)$ the action of $d$ on $\Psi_1(p)$
can be extended to an automorphism of $\Xi^e$. This shows that
$H^u=D^z$ for every $u \in \Xi^e \cap \Gm_2^4(x)$ and it remains to
show that $D^z$ is independent on $z$. Let $v \in \Gm(x)$ be such that
$i_x(v)$ is contained in $\Psi_1(p)$ and put $f=\{x,v\}$. Then by 
(\ref{C} (i)) $\Pi:=\Xi^e \cap \Xi^f$ is a complete graph on 88 vertices. If $R$ is
the stabilizer of $i_x(\Pi \setminus \{x\})$ in $F$, then $R/O_2(R)
\cong L_3(4).S_3$. In view of the remark after (\ref{C}) for an element $r
\in R$ its action on $\Psi_1(p) \cap \Psi_1(q)$ can be extended to an
automorphism of $\Xi^e$ if and only if it can be extended to an
automorphism of $\Xi^f$. Hence $D^z=D^v$ and the result follows by the
connectivity of $\Psi$.\qed

\medskip

\begin{lem} \label{H}
In terms of $(\ref{A})$ let $u_i \in \Gm(x)$ be such that $i_x(u_i)=w_i$
for $i=2,3$ and $4$. Then $u_2,u_3 \in \Gm_2^3(y)$ and $u_4 \in
\Gm_2^4(y)$.\qed 
\end{lem}

\pf For $i=2$ and $4$, in view of (\ref{A}) and (\ref{AB}), the claim
comes as a bi-product of the proof of 
(\ref{G}) for the case $\alpha=3$. By (\ref{A} (iv)) the subgraph
induced by $\Gm(x) \cap \Gm(y) \cap \Gm(u_3)$ has 126 vertices and it
is not a complete graph by (\ref{A} (v)). By (\ref{C} (v)) this can not
happen when $u_3 \in \Gm_2^4(y)$.\qed

\medskip

We need the following result. 

\begin{lem} \label{HI}
Let $^4\widetilde \Xi^e$ be the graph on $2300$ vertices for $^4G \cong
Co_2$ as in $(5.2)$, $f$ be at distance $2$ from $e$, $\Phi$ be the subgraph on $324$
vertices induced by the common neighbours of $e$ and $f$, $\Pi$ be a
complete $22$-vertex subgraph in $\Phi$ and $D=$$^4G(e) \cap $$^4G(f)
\cong U_4(3).2^2$. Then the setwise stabilizer of $\Pi$ in $D$ induces
on $\Pi$ an action of $P\Sigma L_3(4)$.
\end{lem}

\pf Notice that if $\Lm$ is as in (\ref{10}), then $\Phi$ is the subgraph in
the distance 1-or-2 graph of $\Lm$ induced by the orbit $X$ of length 324
of $D$. It is easy to see, using for instance (\ref{16}) that 
the subgraph in $\Lm$ itself induced by $X$ has the following
intersection diagram:

\bigskip

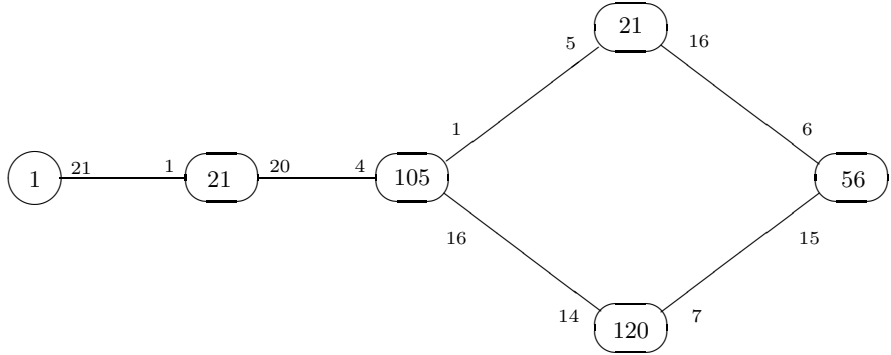
\unitlength 0.80mm
\linethickness{0.4pt}
\begin{picture}(153.66,59.66)
\put(11.99,30.66){\circle{8.00}}
\put(42.99,30.66){\oval(12.00,8.00)[]}
\put(74.66,30.66){\oval(12.00,8.00)[]}
\put(147.66,30.66){\oval(12.00,8.00)[]}
\put(111.00,5.66){\oval(12.00,8.00)[]}
\put(111.00,55.66){\oval(12.00,8.00)[]}
\put(15.99,30.66){\line(1,0){21.00}}
\put(48.99,30.66){\line(1,0){19.67}}
\put(80.00,28.00){\line(4,-3){26.00}}
\put(80.33,33.33){\line(4,3){25.33}}
\put(116.00,52.66){\line(4,-3){26.33}}
\put(11.99,30.33){\makebox(0,0)[cc]{\small 1}}
\put(42.66,30.33){\makebox(0,0)[cc]{\small 21}}
\put(74.66,30.66){\makebox(0,0)[cc]{\small 105}}
\put(111.00,56.00){\makebox(0,0)[cc]{\small 21}}
\put(111.00,5.33){\makebox(0,0)[cc]{\small 120}}
\put(148.00,30.33){\makebox(0,0)[cc]{\small 56}}
\put(19.67,32.33){\makebox(0,0)[cc]{\scriptsize 21}}
\put(34.33,32.66){\makebox(0,0)[cc]{\scriptsize 1}}
\put(52.66,32.66){\makebox(0,0)[cc]{\scriptsize 20}}
\put(65.99,32.66){\makebox(0,0)[cc]{\scriptsize 4}}
\put(82.00,38.66){\makebox(0,0)[cc]{\scriptsize 1}}
\put(101.00,53.00){\makebox(0,0)[cc]{\scriptsize 5}}
\put(122.33,53.33){\makebox(0,0)[cc]{\scriptsize 16}}
\put(140.33,38.66){\makebox(0,0)[cc]{\scriptsize 6}}
\put(82.00,20.66){\makebox(0,0)[cc]{\scriptsize 16}}
\put(100.66,7.66){\makebox(0,0)[cc]{\scriptsize 14}}
\put(122.00,7.66){\makebox(0,0)[cc]{\scriptsize 7}}
\put(140.66,20.67){\makebox(0,0)[cc]{\scriptsize 15}}
\put(116.00,8.33){\line(4,3){26.33}}
\end{picture}

\medskip

By (\ref{12.5} (i)), or otherwise we know that the valency of $\Phi$ is 147
and $147=105+21+21$ is the only way to present the valency as a sum of
suborbits. Now it is a standard fact about dual polar spaces that if
$Y$ is a maximal clique in the distance 1-or-2 graph of $\Lm$, then
either $|Y|=27$ and $Y$ is a quad, or $|Y|=43$ and $Y$ is a vertex
together with its neighbours in $\Lm$. By (\ref{15} (iii)) a quad intersects $X$
in at most 12 vertices, which implies that the clique $\Pi$ 
is a vertex together with its neighbours in $\Lm$ which are contained
in $\Phi$. Hence the stabilizer of $\Pi$ coincides with the stabilizer
in $D$ of a vertex of $\Phi$ and the result follows.\qed

\medskip 

We have proved in (\ref{G}) that for $\alpha \in \{3,4\}$ $H^y$ is
independent on the choice of $y \in \Gm_2^\alpha(x)$ and it remains to
establish the following. 

\begin{lem} \label{I}
There are $y \in \Gm_2^3(x)$ and $u \in \Gm_2^4(x)$ such that $H^y=H^u$.
\end{lem}

\pf By (\ref{H}) we can choose $y \in \Gm_2^3(x)$ and $u \in
\Gm_2^4(x)$ such that $\Pi:=\Sg^y \cap \Sg^u$ contains a complete
subgraph on 22 vertices. It is sufficient to show that 
$S^y[\Pi] \cap S^u[\Pi]$ is not contained in $F^\infty$.
By the proof of (\ref{AB}) we have $F[\Pi] \cong 2 \times 2^{10}.M_{22}.2$,
$S^y[\Pi] \cong 2^{10}.M_{22}.2$ and $F^\infty[\Pi] \cong 2 \times
2^{10}.M_{22}$. So it is sufficient to show that $S^u[\Pi]$ induces in $\Pi$
the natural action of $P\Sigma L_3(4)$. Since $\Pi \subset \Sg^y$, for
any two distinct vertices $p,q$ in $\Pi$ we have $q \in \Psi_2^2(p)$
and hence $\Pi$ maps bijectively onto its image in $\widetilde \Xi^e$,
where $e=\{x,z\}$ is the edge of $\Gm$ such that $u \in \Xi^e$. Hence
the claim follows directly from (\ref{HI}).\qed

\subsection{The second neighbourhood of a vertex}
In this subsection for a given $x \in \Gm$ we analyze the adjacencies
between the vertices in $\Gm_2(x)$. In view of (\ref{3.5}) the results
established in the previous subsection can be summarized as follows.

\begin{prop} \label{J}
Let $x$ be a vertex of $\Gm$, $F_x \cong $$^2E_6(2).S_3$ be the
automorphism group of the subgraph in $\Gm$ induced by $\Gm(x)$. Then there is a
unique subgroup $H_x \cong $$^2E_6(2):2$ of index $3$ in $F_x$ such
that for $\alpha=3$ and $4$ $H_x$ acts transitively on the set 
$$M^\alpha=\{\Gm(x,y) \mid y \in \Gm_2^\alpha(x)\}$$
of $\mu$-graphs as on the cosets of subgroups $Fi_{22}.2$ and
$2^{1+20}.U_4(3).2^2$, respectively. If $X \in M^\alpha$ and $z
\in \Gm(x)$, then $z \in X$ if and only if 
$$H_x(z) \cap H_x[X] \cong 2 \cdot U_6(2):2~~~and~~~H_x(z) \cap H_x[X]
\cong [2^{20}].P\Sigma L_3(4),$$
for $\alpha=3$ and $4$, respectively.\qed
\end{prop}

Notice that for $z \in \Gm(x)$ the action $\pi_z^x$ of $\pi_z$ on
$\Gm(x)$ is an automorphism of the subgraph in $\Gm$ induced by
$\Gm(x)$. Furthermore, $\pi_z^x \in F^\infty$, particularly $\pi_z^x
\in H_x$. 

\medskip

For $\alpha=3$ or $4$ and $X \in M^\alpha$ the subgroup $H_x(z) \cap H_x[X]$ in (\ref{J}) is
determined uniquely up to conjugation in $H_x(z) \cong
2^{1+20}_+:U_6(2).2$. For $X \in M^3$ the subgroup
$(H_x(z) \cap H_x[X])^\prime /Z(H_x(z)) \cong U_6(2)$ is a complement to
$O_2(H_x(z)/Z(H_x(z)))$ in $H_x(z)^\prime /Z(H_x(z))$ and $H_x(z)$ can be
characterized as the stabilizer in $F_x(z)$ of the conjugacy class of
$(H_x(z) \cap H_x[X])^\prime$ in $F_x(z)^\infty$ (compare (\ref{DE})). 
If $e=\{x,z\}$ and $\Xi^e$ is the 4600-vertex subgraph as in
Section~\ref{s5}, then $H_x(z)$ (resp. $H_z(x)$) consists of those
elements of $F_x(z)$ (resp. of $F_z(x)$) whose actions on $\Xi^e \cap
\Gm(x)$ (resp. on $\Xi^e \cap \Gm(z)$) can be extended to automorphisms of
$\Xi^e$. Furthermore, $H_x=H_x(z)F_x^\infty$. This gives an alternative way to
define the subgroups $H_u$ consistently for all $u \in \Gm$. 

\medskip

We need another preliminary result.

\begin{lem} \label{K}
Let $p$ be a vertex of $\Psi$, $\Pi$ 
be the subgraph in $\Psi$ induced by $\Psi_1(p) \cup \Psi_2^2(p)$ 
and $J$ be the automorphism group of
$\Pi$. Then 
$$2^{20}:U_6(2) \cong \bar I < J \le Aut~\bar I \cong 2^{20+2}:U_6(2):S_3$$
where $\bar I$ is as in $(\ref{DE})$, so that $Out~\bar I \cong S_4$ acts
naturally on the classes of complements in $\bar I$ to $O_2(\bar I)$.
\end{lem} 

\pf By (\ref{1}) the stabilizer of $p$ in $F$ induces on $\Pi$ an
action of $2^{20}:U_6(2).S_3$, so to prove the lemma it is sufficient
to show that $J$ is contained in $Aut~\bar I$. 
Let $\Upsilon$ be a symplecton containing $p$ and $A$ be the action
induced by $J[\Upsilon]$ on $\Upsilon$. 
Then $A$ is contained in the stabilizer of $p$ in 
$Aut~\Upsilon \cong \Omega_8^-(2).2$ and $A$ is contained in the action
induced on $\Upsilon$ by its stabilizer in $F(p)$. 
Hence $A \cong 2^6.\Omega_6^-(2).2$ and $A$ acts faithfully on
$\Upsilon \cap \Psi_1(p)$. By (\ref{3}) $\Upsilon$ is uniquely determined by
its intersection with $\Psi_1(p)$ and hence $J$ acts faithfully on
$\Psi_1(p)$. 
Let $V$ be the kernel of the action of $J$ on the set $L$
of 891 lines containing $p$ and $\widetilde J=J/V$. In every symplecton
containing $p$ there are exactly 27 lines containing $p$. 
The collection of 693
such 27-element subsets of $L$ defines on $L$ the structure $\cD$ of
dual polar space of $U_6(2)$. This structure is preserved by $J$ and
hence $\widetilde J \le U_6(2).S_3$. Since every line in $L$ contains
besides $p$ exactly two points, the kernel $V$ is an elementary abelian
2-group. On every line from $L$ the group $V$ induces an action of order at most 2 and on
the 27 lines in a symplecton $\Upsilon$ it induces an action of order at most
$2^6=|O_2(A)|$. Hence on a plane containing $p$ the subgroup
$V$ induces an action of order at most $2^2$. This shows that the dual
of $V$ is a representation module of $\cD$ and the latter was proved in
\cite{Y94} to be of order $2^{22}$. Hence $J^\infty \cong \bar I$ and
the result follows.\qed

\medskip

For $z \in \Gm(x)$ put $\Omega=(\Gm(x) \cap \Gm(z))\setminus \{x,z\}$. The mappings $i_x$
and $i_z$ restricted to $\Omega$ are bijections onto the vertex set of
the graph $\Pi$ as in (\ref{K}). 
If $u \in \Omega$, then by (\ref{6}) the type of 
$\{x,z,u\}$ is well defined and we can put
$$\Omega_1=\{u \in \Omega \mid \{x,z,u\} {\rm ~is~short}\}$$
$$\Omega_2^2=\{u \in \Omega \mid \{x,z,u\} {\rm ~is~long}\}.$$
If $e=\{u,v\}$ is an edge of $\Gm$ contained in $\Omega$, then the
triangles $\{x,u,v\}$ and $\{z,u,v\}$ might or might not be of the same
type. Define $(\tau_x,\tau_z)$ to be the type of $e$, where
$$\tau_x=s {\rm ~or~} l {\rm ~if~} \{x,u,v\} {\rm ~is~} short{\rm ~or~}
long{\rm,~respectively}$$
and $\tau_z$ is defined similarly. By (\ref{9.75}) we have the
following. 

\begin{lem} \label{M}
Let $e=\{u,v\}$ be an edge of $\Gm$ contained in $\Omega$ and
$(\tau_x,\tau_z)$ be the type of $e$. Then
\begin{itemize}
\item[{\rm (i)}] if $e \subset \Omega_1$ or $e \subset \Omega_2^2$, then
$\tau_x=\tau_z$; 
\item[{\rm (ii)}] if $u \in \Omega_1$ and $v \in \Omega_2^2$, then
$\tau_x \ne \tau_z$.\qed
\end{itemize}
\end{lem}

\begin{lem} \label{N}
Let $K$ be the group of automorphisms of the subgraph in $\Gm$ induced
by $\Omega$, preserving the above defined types $(\tau_x,\tau_z)$ of
edges. Then $K \cong J =Aut~\Pi$, where $\Pi$ is as in $(\ref{K})$.
\end{lem}

\pf For $\gamma=x$ or $z$ let $\Omega^\gamma$ be the graph on $\Omega$,
whose edges are the edges of $\Gm$ contained in $\Omega$ for which 
$\tau_\gamma=s$. Then $i_\gamma$ induces an isomorphism of
$\Omega^\gamma$ onto the graph $\Pi$. In view of (\ref{M}) every
automorphism of $\Pi$ can be realized as a type-preserving automorphism
of the subgraph induced by $\Omega$. Hence the result.\qed 

\medskip

Let $^x\cY$ be the set of $x$-equivalence classes of vertices in
$\Gm_2(x)$. By (\ref{J}) we can define the action of $H_x$ on $^x\cY$ by
the following rule: if $\{y,y^\prime\} \in $$^x\cY$, with
$\Gm(x,y)=\Gm(x,y^\prime)=X$, then $\{y_1,y_1^\prime\}=\{y,y^\prime\}^h$ 
for $h \in H_x$ if and only if 
$\Gm(x,y_1)=\Gm(x,y_1')=X^h$. By (\ref{F}) and (\ref{G}) this
action is well defined. For $z \in \Gm(x)$ let $H_x(z)$ be the
stabilizer of $z$ in $H_x$ and $^x\cY_z$ be the set of classes from $^x\cY$
contained in $\Gm(z)$. Since $\{y,y^\prime\} \in $$^x\cY_z$ whenever $z
\in \Gm(x,y)$, $H_x(z)$ stabilizes $^x\cY_z$ as a whole and we can
consider the action of $H_x(z)$ on $^x\cY_z$. 

\begin{lem} \label{L}
The above defined action of $H_x(z)$ on $^x\cY_z$ coincides with
the natural action of $H_z(x)$ on $^x\cY_z$.
\end{lem}

\pf With $\Omega=(\Gm(x) \cap \Gm(z)) \setminus \{x,z\}$ as above, 
let $\{y,y^\prime\} \in $$^x\cY_z$ and $X=\Gm(x,y)$. By
(\ref{EF}) and (\ref{8.5}) the stabilizer of $\{y,y^\prime\}$ in
$H_z(x)$ coincides with the setwise stabilize of $\Gm(y) \cap \Omega$
and the stabilizer of $X$ in $H_x(z)$ coincides with the stabilizer of
$X \cap \Omega$. Furthermore $\Gm(y) \cap \Omega=X \cap \Omega$, which 
implies that the actions of $H_x(z)$ and $H_z(x)$ on
$^x\cY_z$ are determined by their natural actions on $\Omega$. Thus all
we have to show is that the action $A_x$ of $H_x(z)$ on $\Omega$
coincides with the action $A_z$ of $H_z(x)$ on this set. Notice that
$A_x \cong A_z \cong 2^{20}:U_6(2).2$ and both $A_x$ and $A_z$ are
subgroups in the group $K$ of type-preserving automorphisms of the
subgraph in $\Gm$ induced by $\Omega$ (compare (\ref{N})). Let $t_0 \in
\Gm(x)$ with $i_x(t_0) \in \Psi_3(i_x(z))$, $t_1 \in \Gm(z)$ with
$i_z(t_1) \in \Psi_3(i_z(x))$. For $i=0$ and $1$ let $U_i$ be the
stabilizer of $\Gm(t_i) \cap \Omega$ in $A_x$. Then $U_0 \cong U_1 
\cong U_6(2).2$ and $U_0$ and $U_1$ belong to different classes of
complements to $O_2(A_x)$. By the paragraph after (\ref{J}) $A_x$ 
is characterized as the stabilizer in $J \cong 2^{20}:U_6(2).S_4$ of the classes in $J^\infty$
containing $U_0^\prime$ and $U_1^\prime$. Since $A_z$ is isomorphic to
$A_x$ and possesses the same properties, the result follows.\qed 
 
\medskip

Now we can establish some further properties of the permutation
$\pi_x$. 

\begin{lem} \label{pi2}
If $s,w \in \Gm(x)$ and $s \in \Gm(w)$ then the restriction of $\pi_x$
to the set $\Gm(x) \cup \Gm(s) \cup \Gm(w)$ is an automorphism of the
subgraph in $\Gm$ induced by this set.
\end{lem}

\pf In terms of (\ref{L}) let $d \in \Gm(x) \cap \Gm(z)$. Then the
action of $\pi_d$ on $\Gm(x)$ is an element in $H_x(z)$ while its
action on $\Gm(z)$ is an element of $H_z(x)$. Hence the result follows
from (\ref{L}).\qed

\medskip

By now with every vertex $z \in \Gm$ we have associated various actions
on $\Gm$. On the one hand in (\ref{defi}) we have defined the
permutation $\pi_z$ of the vertex-set of $\Gm$. On the other hand, if
$x \in \Gm(z)$, then the action $\pi_z^x$ of $\pi_z$ on the set
$\Gm(x)$ is an element of $H_x$ and an
action of $\pi_z^x$ on the set $^x\cY$ of pairs of $x$-equivalent
vertices in $\Gm_2(x)$ is defined in the paragraph before (\ref{L}).
The natural question is whether these actions are consistent. A partial
affirmative answer is given in the following lemma.

\begin{lem} \label{part}
Let $x \in \Gm$, $z \in \Gm(x)$, $y \in \Gm_2(x)$ and suppose that $z$
is adjacent in $\Gm$ to a vertex $u \in \Gm(x,y)$. Then
$\Gm(y^{\pi_z},x)$ is the image of $\Gm(x,y)$ under $\pi_z^x$, i.e.
$\Gm(x,y)^{\pi_z^x}=\Gm(y^{\pi_z},x)$.
\end{lem}

\pf Under the hypothesis $x,y,z \in \Gm(u)$. Hence $y^{\pi_z}$
coincides with the image of $y$ under $\pi_z^u$ and $\Gm(y^{\pi_z},x)
\cap \Gm(u)$ is the image of $\Gm(y,x) \cap \Gm(u)$ under $\pi_z^x$
(equivalently under $\pi_z^u$). By (\ref{EF}) $\Gm(y^{\pi_z},x)$ is
uniquely determined by $\Gm(y^{\pi_z},x) \cap \Gm(u)$ and hence the
result.\qed 

\subsection{$\pi_x$ is an automorphism of $\Gm$}
In this subsection we prove Proposition~\ref{ult}. We need to show that
whenever $\{y,u\}$ is an edge of $\Gm$, its image under $\pi_x$ is
again an edge of $\Gm$. Since $\pi_x$ fixes every vertex in $\Gm
\setminus \Gm_2(x)$, we can assume that at least one of the vertices on
the edge, say $y$, is in $\Gm_2(x)$. Then $y^{\pi_x}=y^\prime$ is
$x$-equivalent to $y$ in $\Gm_2(x)$. If there is $z \in \Gm(x)$ such
that $\{y,u\} \subset \{z\} \cup \Gm(z)$, then
$\{y^{\pi_x},u^{\pi_x}\}$ is an edge by (\ref{pi0}). Hence we assume
that 
$$\Gm(x) \cap \Gm(y) \cap \Gm(u)=\emptyset.$$
By (\ref{B} (ii)) this implies that $y \not \in \Gm_2^3(x)$ and so $y
\in \Gm_2^4(x)$. Furthermore by the proof of (\ref{E}) for $\alpha=4$
we conclude that $x \in \Xi^e$ for the edge $e=\{y,y^\prime\}$. 

>From now on we follow notation introduced in the paragraph after the
proof of (\ref{13.5}) (with the roles of $x$ and $y$ interchanged). As
usual $i_y$ is a fixed isomorphism of the subgraph in $\Gm$ induced by
$\Gm(y)$ onto the graph $\Dl$. Let
$$\Omega=\{i_y(v) \mid v \in \Gm(x,y) \},$$
$I \cong 2^{1+20}_+:U_4(3).2^2$ be the setwise stabilizer of $\Omega$
in $F=Aut~\Dl$ and $\Omega_0,\Omega_1=\Omega,...,\Omega_7$ be the
orbits of $I$ on the vertex-set of $\Dl$, where $\Omega_0=\{p\}$,
$p=i_y(y^\prime)$, $\Omega_1 \cup \Omega_2=\Psi_1(p)$, $\Omega_3 \cup
\Omega_4=\Psi_2^2(p)$, $\Omega_5 \cup \Omega_6=\Psi_2^4(p)$,
$\Omega_7=\Psi_3(p)$ (the notation is as in (\ref{16})). 

If $i_y(u) \in \Omega_\beta$ for $\beta \in\{0,1,2,4,5,6\}$, then by
(\ref{17}) $i_y(u)$ is adjacent in $\Dl$ to a vertex from $\Omega$
which means that $u$ is adjacent in $\Gm$ to a vertex from $\Gm(x,y)$
and this is not the case we consider. Hence we assume that $i_y(u) \in
\Omega_3 \cup \Omega_7$. For $\beta=3$ or $7$ put $\Phi^x_\beta(y)=\{w
\mid w \in \Gm(y), i_y(w) \in \Omega_\beta\}$.

\begin{lem} \label{prel}
The following assertions hold:
\begin{itemize}
\item[{\rm (i)}] if $q \in \Omega_3$ then $\Psi_3(q) \cap \Omega=\emptyset$;
\item[{\rm (ii)}] if $q \in \Omega_7$ then $|\Psi_3(q) \cap
\Omega|=\frac{1}{2}|\Omega|$.
\end{itemize}
\end{lem}

\pf If $q \in \Omega_3$ and $r \in \Omega$, then by (\ref{1} (ii),
(iii)) $r \in O_2(F(p))$ and $q \in F(p)$ which means that $\la q,r\ra$
has even order. If $q \in \Omega_7$ then the claim is by (\ref{1}
(v)).\qed 

\begin{lem} \label{lem}
Suppose that $u \in \Phi^x_7(y)$. Then $u \in \Gm_2^4(x)$ and $y \in
\Phi^x_7(u)$. 
\end{lem}

\pf By (\ref{prel}) there is $w \in \Gm(x,y)$, such that $u \in
\Gm_2^3(w)$. Then $x \in \Gm(w) \setminus \Gm(w,u)$ and by (\ref{B}
(ii)) $u \in \Gm_2(x)$. Since $\Gm(x) \cap \Gm(y) \cap
\Gm(u)=\emptyset$, we have $u \in \Gm_2^4(x)$. By (\ref{B} (ii)) $x$ is
adjacent to a vertex $z \in \Gm(u,w)$. Then both  $y$ and $z$ are
contained in $\Gm(u,w)$ and are at distance 2 in $\Gm$. By (\ref{13.5})
this means that $y \in \Gm_2^3(z)$ and by (\ref{prel}) $y \in
\Phi^x_7(u)$.\qed 

\medskip

Now we can easily handle the case $u \in \Phi^x_3(y)$.

\begin{lem} \label{hand}
Suppose that $u \in \Phi^x_3(y)$. Then the image of $\{y,u\}$ under
$\pi_x$ is an edge of $\Gm$.
\end{lem}

\pf The crucial observation is that in the considered situation
$i_y(u)$ is adjacent in $\Dl$ to $p=i_y(y^\prime)$, which means that
$u$ is adjacent to $y^\prime$ in $\Gm$. If $u \in \Gm_3(x)$ (which is
the case in the Baby Monster graph), then $u^{\pi_x}=u$ and hence
$\{y^\prime,u\}$ is an edge by the above observation. Suppose that $u
\in \Gm_2(x)$. Then $u \in \Gm_2^4(x)$ by (\ref{B} (ii)), since $\Gm(x)
\cap \Gm(y) \cap \Gm(u)=\emptyset$. Furthermore by (\ref{17}) and (\ref{cor}) $y \in
\Phi^x_3(u)$ or $y \in \Phi^x_7(u)$ and the latter case is impossible by
(\ref{lem}). Hence $y \in \Phi^x_3(u)$ and by the obvious symmetry the
subgraph in $\Gm$ induced by $\{y,y^\prime,u,u^\prime\}$ (where $u^\prime$
is the $x$-equivalent to $u$) is complete. Hence the result.\qed

\medskip

Now we turn to the final case, when $u \in \Phi^x_7(y)$. We follow
notation from the proof of (\ref{lem}), so that $w \in \Gm(x,y)$, $u
\in \Gm_2^3(w)$, $z \in \Gm(x) \cap \Gm(w,u)$. Then $x,y,z \in \Gm(w)$
and $u \in \Gm_2(w)$. Hence $\pi_x^w$ is an element of $H_w$, $y \in
\Gm(u,w)$ and $x$ is
adjacent to $z \in \Gm(u,w)$. Then by (\ref{part})
$\Gm(u^{\pi_x},w)$ is the image of $\Gm(u,w)$ under $\pi_x^w$, in
particular it contains $y^{\pi_x}$. Hence $\{y^{\pi_x},u^{\pi_x}\}$ is
an edge of $\Gm$. This completes the proof of (\ref{ult}).

\subsection{Proof of Theorems~\ref{aaa} and \ref{bbb}}
By Proposition~\ref{ult} for every $x \in \Gm$ the permution $\pi_x$ is
an automorphism of $\Gm$. The actions $\pi_z^x$ of the automorphisms
$\pi_z$ on $\Gm(x)$ taken for all $z \in \Gm(x)$ generate $F_x^\infty
\cong $$^2E_6(2)$, hence Theorem~\ref{aaa} follows from Proposition~\ref{a}. Now
Theorem~\ref{bbb} follows from (\ref{5}) and Theorem~\ref{aaa}.

\newpage

\begin{center}

{\Large Acknowledgments}

\end{center}

\bigskip

The authors are sincerely thankful to Corinna Wiedorn, who have read
very carefully the previous version of the paper and made a number of
very helpful suggestions on how to improve the exposition.

\bigskip

A.A.~Ivanov,

Department of Mathematics,

Imperial College,

180 Queen's Gate,

London, SW7 2BZ, UK

\medskip

D.V.~Pasechnik, 

SSOR/TWI,  

Delft University of Technology,

Mekelweg 4,

2628 CD Delft,

The Netherlands

\medskip

S.V.~Shpectorov,

Department of Mathematics and 

Statistics,

Bowling Green State University,

Bowling Green, OH 43403, USA


\bigskip



\begin{thebibliography}{WWWWWW}

\bibitem[BCN89]{BCN89} A.E.~Brouwer, A.M.~Cohen and A.~Neumaier, {\it
Distance-Regular Graphs}, Springer-Verlag, Berlin, 1989.

\bibitem[B85]{B85} F.~Buekenhout, Diagram geometries for sporadic
groups, In: {\it Finite Groups---Coming of Age}, J. McKay ed., {\it
Contemp. Math.} {\bf 45} (1985), 1--32.

\bibitem[Coh83]{Coh83} A.M.~Cohen, Points and lines in metasymplectic
spaces, {\it Ann. Discr. Math.} {\bf 18} (1983), 193--196.

\bibitem[CC88]{CC88} A.M.~Cohen and B.~Cooperstein, The $2$-spaces of
the standard $E_6(q)$-module, {\it Geom. Dedic.} {\bf 25} (1988),
467--480.

\bibitem[Coo83]{Coo83} B.N.~Cooperstein, The geometry of root subgroups in
exceptional groups II, {\it Geom. Dedic.} {\bf 15} (1983), 1--45.

\bibitem[ATLAS]{ATLAS} J.H.~Conway, R.T.~Curtis, S.P.~Norton,
R.A.~Parker and R.A.~Wilson, {\it Atlas of Finite Groups}, Clarendon
Press, Oxford, 1985.

\bibitem[C94]{C94} H.~Cuypers, A graphical characterization of $Co_2$,
Preprint 1994.

\bibitem[C99]{C99} H.~Cuypers, Extended near hexagons and line system,
Preprint, 1999.

\bibitem[CP92]{CP92} H.~Cuypers and A.~Pasini, Locally polar geometries
with affine planes, {\it Europ. J. Combin.} {\bf 13} (1992), 39--57.

\bibitem[HS00]{HS00} J.I.~Hall and S.V.~Shpectorov, Rank 3 P-geometries,
{\it Geom. Dedic.}, {\bf 82} (2000), 139--169.

\bibitem[I92]{I92} A.A.~Ivanov, A geometric characterization of
Fischer's Baby Monster, {\it J. Algebraic Comb.} {\bf 1} (1992),
43--65.

\bibitem[I94]{I94} A.A.~Ivanov, Presenting the Baby Monster, {\it J.
Algebra} {\bf 163} (1994), 88--108.

\bibitem[I95]{I95} A.A.~Ivanov, On geometries of Fischer groups, {\it
Europ. J. Combin.} {\bf 16} (1995), 163--183.

\bibitem[I99]{I99} A.A.~Ivanov, {\it Geometry of Sporadic Groups I.
Petersen and Tilde Geometries}, Cambridge Univ. Press, Cambridge, 1999.

\bibitem[ILLSS]{ILLSS} A.A.~Ivanov, S.A.~Linton, K.~Lux, J.~Saxl and
L.H.~Soicher, Distance-transitive representations of the sporadic
groups, {\it Comm. Algebra} {\bf 23} (1995), 3397--3427.

\bibitem[IS96]{IS96} A.A.~Ivanov and J.~Saxl, The character table of
$^2E_6(2)$ acting on the cosets of $Fi_{22}$, In: {\it Advanced Studies
in Pure Math.} {\bf 24} (1996), 165--196.

\bibitem[L93]{L93} R.~Lawther, Folding actions, {\it Bull. London Math.
Soc.} {\bf 25} (1993), 132--144.

\bibitem[LS98]{LS98} M.W.~Liebeck and G.M.~Seitz, On the subgroup
structure of exceptional groups of Lie type, {\it Trans. Amer. Math.
Soc.} {\bf 350} (1998), 3409--3482.

\bibitem[N92]{N92} S.P.~Norton, Constructing the Monster, In: {\it
Groups, Combinatorics and Geometry}, Durham 1990, M.~Liebeck and
J.~Saxl eds., London Math. Soc. Lect. Notes, {\bf 165}, Cambridge Univ.
Press, Press, 1992, pp. 63--76.

\bibitem[Pase94]{Pase94} D.V.~Pasechnik, Geometric characterization of the
sporadic groups $Fi_{22}$, $Fi_{23}$ and $Fi_{24}$, {\it J. Combin.
Theory} (A) {\bf 68} (1994), 100--114.

\bibitem[Pase95]{Pase95} D.V.~Pasechnik, Extended polar spaces of rank at
least $3$, {\it J. Combin. Theory} (A) {\bf 72} (1995), 232--242.

\bibitem[Pasi94]{Pasi94} A.~Pasini, {\it Diagram Geometries}, Clarendon
Press, Oxford, 1994.

\bibitem[PS97]{PS97} C.E.~Praeger and L.H.~Soicher, {\it Low Rank
Representations and Graphs for Sporadic Groups}, Cambridge Univ. Press,
Cambridge, 1997.

\bibitem[T74]{T74} J.~Tits, {\it Buildings of Spherical Type and Finite
$BN$-pairs}, Lect. Notes Math., {\bf 386}, Springer-Verlag, Berlin
1974. 

\bibitem[Y94]{Y94} S.~Yoshiara, On some extended dual polar spaces I.,
{\it Europ. J. Combin.} {\bf 15} (1994), 73--86. 

\end{thebibliography}
\end{document}